\documentclass[a4paper, 11pt]{article}

\usepackage{amssymb}
\usepackage{mathrsfs}
\usepackage{amsfonts}
\usepackage{graphicx,graphics,epic}

\title{\textbf{Ternary algebras with braided statistics}}
\date{}
\author{
 Azzouz.ZINOUN\\
 To my father's memory.\\
{Laboratoire Phlam,UMR-CNRS 8523 UFR de physique,}\\{\small Cit\'e
Scientifique~Universit\'{e}
des sciences et Technologies de Lille,}\\{\small 59655 Villeneuve d'Ascq cedex France.}\\
{\small E-mail : Azzouz.Zinoun@univ-lille1.fr}}
\begin{document}
\maketitle
\begin{abstract}
Algebraic relations that characterize quantum statistics (
Bose-Einstein statistic, Fermi-Dirac statistic, supersymmetry,
parastatistic, anyonic statistic,... )are reformulated herein in
terms of a new algebraic structure, which we call para-algebra .
\end{abstract}
\section{Introduction}
There have been suggestions for quantization procedures that lead to
particle statistics different from the well known Bose-Einstein and
Fermi-Dirac Statistics. Elaborate mathematical developpement are
needed to justify such quantization procedures. An interesting
example is parastatistics which was introduced by H.S~Green
\cite{gr} in 1953. It is well known that the idea of supersymmetry
\cite{ne} concerns bosons and fermions and the mathematic that
describe this field is Lie superalgebras\cite{Ka}, or equivalently
$\mathbb {Z}$$_2$-graded algebras, which is based on binary
relations. In an earlier paper \cite{Zi}, we have shown that there are general
mathematical structures that
represent different statistics, and we have established the connection between them. In practice, braided groups are simply an analog of
supergroups with $\pm1$ Bose-Fermi statistics replaced by braid
statistics. There are indications that particles of braid
statistics arise in low dimension quantum field theory.\\
On the other hand the parasupersymmetric \cite{Oh} quantum
mechanics, is related to parabosons and parafermions and the
mathematical theory need to study this field is based on ternary
operations which we will call para-algebra.\\
In this paper, we attempt to construct the theory of para-algebras,
this terminology comes from physics in the context of
parastatistics. Parasuperalgebras or equivalently $\mathbb
Z_{2}$-graded para-algebras are para-algebras where the bilinear map
$\sigma$ defined in paper \cite{Zi} takes its values in $\mathbb
Z_{2}$; from this formalism we derive the algebraic structures which
lead to quantum statistics, specially to generalized parastatistics
\cite{Bi}, which we will call Lie parasuperalgebra.  In section
$\bf{2}$ we give the basic definitions concerning braided tensor
algebra \cite{Ma, Ca, PC}.  In section $\bf{3}$  we  illustrate our
method, by defining maps on a braided tensor algebra,then we
construct ideals and envelopping algebras to define ternary maps
which give the para-algebras.  In section  $\bf{4}$ , we derive the
formalism of the generalized parastatistics \cite{Bi}. In section 5, we show that these ternary maps are related to Schur fonctor or Weyl module.

\section{Preliminaries}
\subsection{definitions}
Let $\Gamma$ be an n dimensional vector space over a field $K$ of
characteristic $0$, and let $(A,\otimes,I)$ be a tensorial category
$\cite{Ma}$,\cite{Ca},\cite{PC}; the objects
 of $A$ are denoted by $A_{v_{i}}, i = 0, 1,\ldots$; $A_{v_{0}}=I$; $ v_{i}\in \Gamma $ and $I$
 is the unit element for the operation $\otimes$, for simplicity
 $A_{v_{i}}\otimes I=I \otimes A_{v_{i}} = A_{v_{i}}$ ( this notation will be clear
later; subsequently the $A_{v_{i}}$ will denote vector spaces).\\
$(A,\otimes,I)$ is a braided tensorial category if there exist
natural isomorphisms called braiding defined between any two objects
of $A$ such that
\begin{eqnarray}
\Psi_{A_{v_{i}},A_{v_{j}}} : A_{v_{i}} \otimes A_{v_{j}} \rightarrow
A_{v_{j}} \otimes  A_{v_{i}}
\end{eqnarray}
Let $\left\{e_{i}\right\}$ and $\left\{f_{j}\right\}$ be,
respectively, the bases of the vector spaces $A_{v_{i}}$ and
$A_{v_{j}}$. Then :
\begin{eqnarray}
\Psi_{A_{v_{i}},A_{v_{j}}}( e_{i} \otimes f_{j} )& =& \Psi^{m n}_{i
j} f_{m} \otimes e_{n}
\end{eqnarray}
(The summation on the repeated indices is to be
understood).Therefore
\begin{eqnarray}
\Psi_{A_{v_{i},I}} = \Psi_{I,A_{v_{i}}}
 = id_{A_{v_{i}}}
\end{eqnarray}
\begin{eqnarray}
\Psi_{A_{v_{i}},A_{v_{j}}\otimes A_{v_{k}}} = \left(id_{A_{v_{j}}}
\otimes \Psi_{A_{v_{i}},A_{v_{k}}}\right)
\circ\left(\Psi_{A_{v_{i}},A_{v_{j}}} \otimes id_{A_{v_{k}}}\right)
\end{eqnarray}
\begin{eqnarray}
\Psi_{A_{v_{i}}\otimes A_{v_{j}}, A_{v_{k}}} =
\left(\Psi_{A_{v_{i}},A_{v_{k}}} \otimes
id_{A_{v_{j}}}\right)\circ\left(id_{A_{v_{i}}}\otimes
\Psi_{A_{v_{j}},A_{v_{k}}}\right)
\end{eqnarray}
\begin{eqnarray}
\label{e6} \Psi^{m n}_{i j}\Psi^{r q}_{m n} =
\delta^{r}_{i}\delta^{q}_{j}~~ unitarity~ condition
\end{eqnarray}
The property $(4)$ is equivalent to the following triangular diagram :\\
\vspace{0.5cm}\\
\unitlength=4mm \put(12,1){$\Psi_{A_{v_{i}},A_{v_{j}}\otimes
A_{v_{k}}}$} $A_{v_{i}} \otimes A_{v_{j}}\otimes A_{v_{k}}$
\put(0,0){{\vector(4,0){15}} $A_{v_{j}} \otimes A_{v_{k}}\otimes
A_{v_{i}}$} \put(-3.5,-1){\vector(3,-1){9}}\put(3.1,-5){$A_{v_{j}}
\otimes A_{v_{k}}\otimes A_{v_{i}}$}\put(8,-4){\vector(3,1){9}}
\put(-2.2,-2.9){$\Psi_{A_{v_{i}},A_{v_{j}}}\otimes id_{A_{v_{k}}}$}
\put(10,-2.9){$id_{A_{v_{j}}} \otimes \Psi_{A_{v_{i}},A_{v_{k}}}$}
\vspace{0.5cm}\\
The property $(5)$ is equivalent to the following triangular
diagram:
\vspace{0.5cm}\\
\unitlength=4mm \put(12,1){$\Psi_{A_{v_{i}}\otimes A_{v_{j}},
A_{v_{k}}}$} $A_{v_{i}} \otimes A_{v_{j}}\otimes A_{v_{k}}$
\put(0,0){{\vector(4,0){15}} $A_{v_{k}} \otimes A_{v_{i}}\otimes
A_{v_{j}}$} \put(-3.5,-1){\vector(3,-1){9}}\put(3.1,-5){$A_{v_{i}}
\otimes A_{v_{k}}\otimes A_{v_{j}}$}\put(8,-4){\vector(3,1){9}}
\put(-2.2,-2.9){$id_{A_{v_{k}}}\otimes\Psi_{A_{v_{j}},A_{v_{k}}} $}
\put(10,-2.9){$ \Psi_{A_{v_{i}},A_{v_{k}}}\otimes id_{A_{v_{j}}}$}
\vspace{0.5cm}\\
In a braided tensorial category, we have the following identity:

\begin{eqnarray}
\left(\Psi_{A_{v_{j}},A_{v_{k}}}\otimes
id_{A_{v_{i}}}\right)\left(id_{A_{v_{j}}} \otimes
\Psi_{A_{v_{i}},A_{v_{k}}}\right)\left(\Psi_{A_{v_{i}},A_{v_{j}}
\otimes A_{v_{k}}}\right)~\nonumber\\=~\left(id_{A_{v_{k}}}\otimes
\Psi_{A_{v_{i}},A_{v_{j}}}\right)
\left(\Psi_{A_{v_{i}},A_{v_{k}}\otimes
A_{v_{j}}}\right)\left(id_{A_{v_{i}}}\otimes
\Psi_{A_{v_{j}},A_{v_{k}}}\right)
\end{eqnarray}
which is the generalized Yang-Baxter identity : \vskip 1cm
\begin{figure}[!h]
\begin{center}
\includegraphics[width=10cm,height=8.5cm]{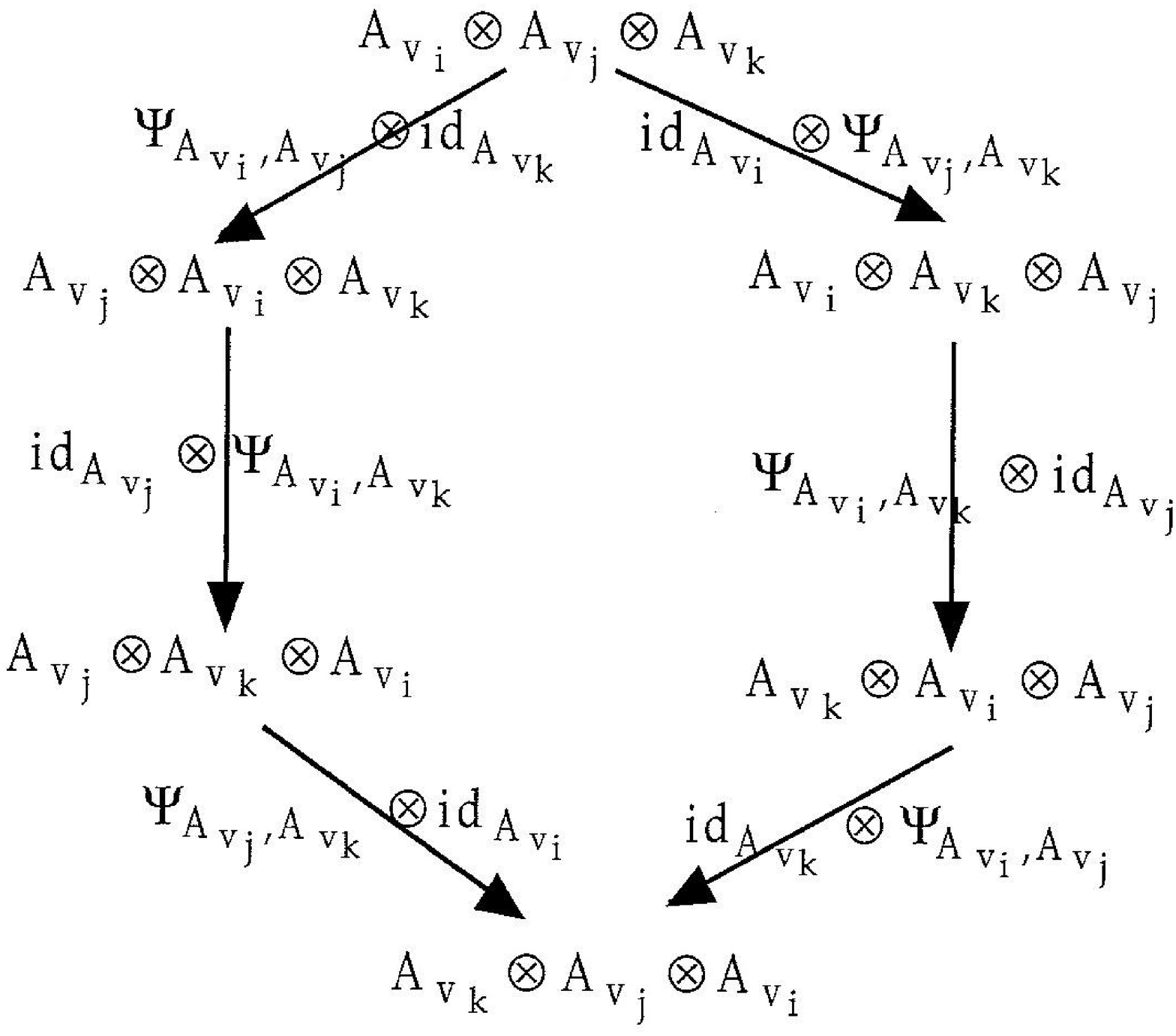}
%\includegraphics[width=12cm]{figure=1.eps}
%\subfigure{\epsfig{figure=1.eps,width=10cm,height=8.5cm}}
\end{center}
\end{figure}

\section{Ternary maps}
\subsection{Definitions}
Let $A$ be a $\Gamma$-graded vector space \cite{Zi} $A=\oplus
A_{v_{i}}$; $v_{i}\in \Gamma$, $i=1\ldots n$, $A_{v_{i}}$ is a
vector subspace of $A$. Let $T(A)=\oplus_{P\geq 0} A^{\otimes P}$ be
the tensor algebra constructed from the vector space $A$; since $A$
is braided, $T(A)$ is also braided. From the properties a) and b) in
section 2 one can see that there are only two ways for braiding a
3-fold tensor product. It is useful to write somme notations; this
consists of writing  all isomorphisms pointing downwards with
$\Psi=$ \unitlength=4mm \put(0,1){\vector(4,-3){2,4} $$}
\put(2.5,1){\line(-4,-3){1,1}$$}
\put(1.15,-0.1){\vector(-4,-3){1}$$}~~~~~~~~~,
$\Psi^{-1}=$\unitlength=4mm \put(0,1){\line(4,-3){1,1} $$}
\put(2.5,1){\vector(-4,-3){2,4}$$}
\put(1.3,-0.1){\vector(4,-3){1}$$}~~~~~~~~~as braids.

%$$ \Psi: A\otimes A\otimes A \rightarrow A\otimes A\otimes A $$
%
%%%%%%%%%%%%%%%%%%%%%%%%%%%%%%%%%%%%%%%%%%%%%%%%%%%%%%
%%%%%%%%%%%%%%%%%%%%%%%%%%%%%%%%%%%%%%%%%%%%%%%%%%%%%%%
%%%%%%%%%%%%%%%%%%%%%%%%%%%%%%%%%%%%%%%%%%%%%%%%%%%%%%%%%%
%%%%%%%%%%%%%%%%%%%%%%%%%%%%%%%%%%%%%%%%%%%%%%%%%%%%%%%%%%%
We need the following maps:\\
\subsection{Left ternary maps}
 By braiding from the left we construct the following maps \vspace{0.4cm}
\begin{figure}[!ht]
{\begin{minipage}{0.45\linewidth} \centering
\includegraphics[width=6cm,height=2.7cm]{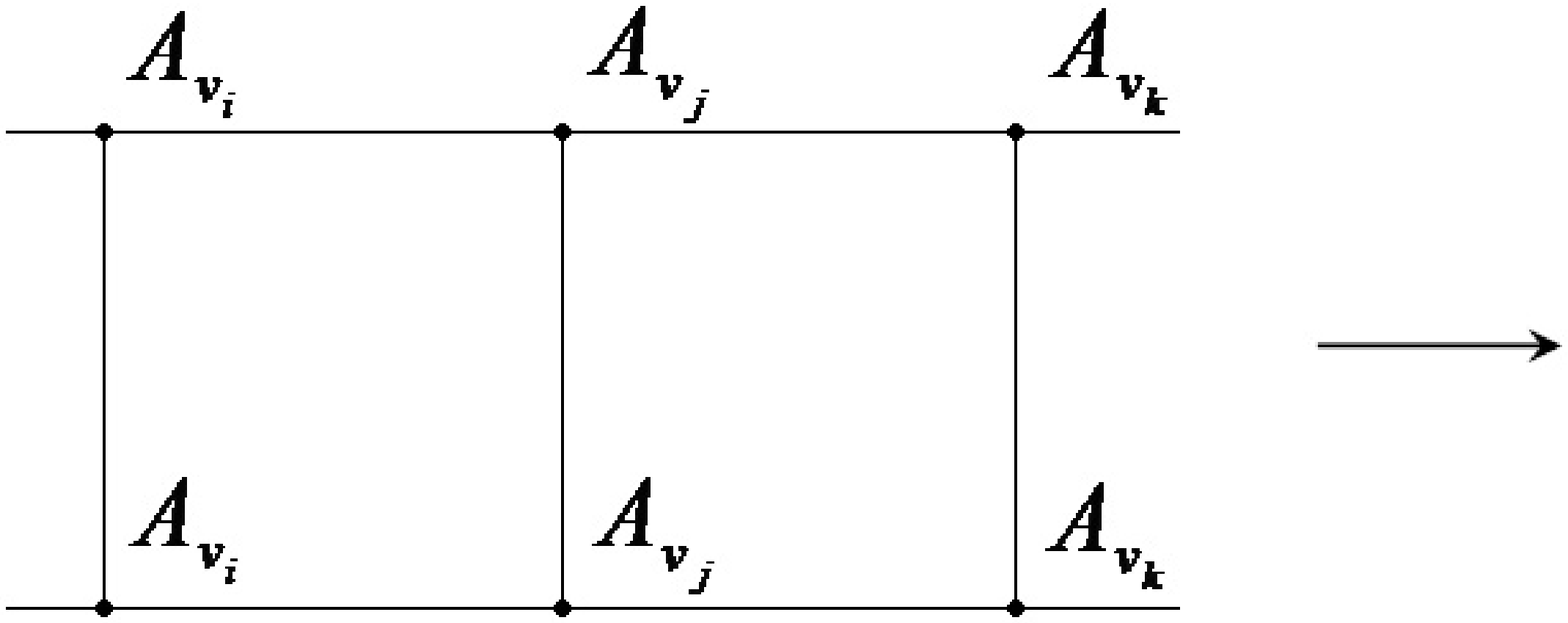}
\end{minipage}\hfill}
{\begin{minipage}{0.5\linewidth}$id=id_{A_{v_{i}}}\otimes
id_{A_{v_{j}}}\otimes id_{A_{v_{k}}}$
\end{minipage}\hfill}
%\end{figure}
\vskip 0.5cm
%\begin{figure}[!ht]
{\begin{minipage}{0.45\linewidth} \centering
\includegraphics[width=6cm,height=2.7cm]{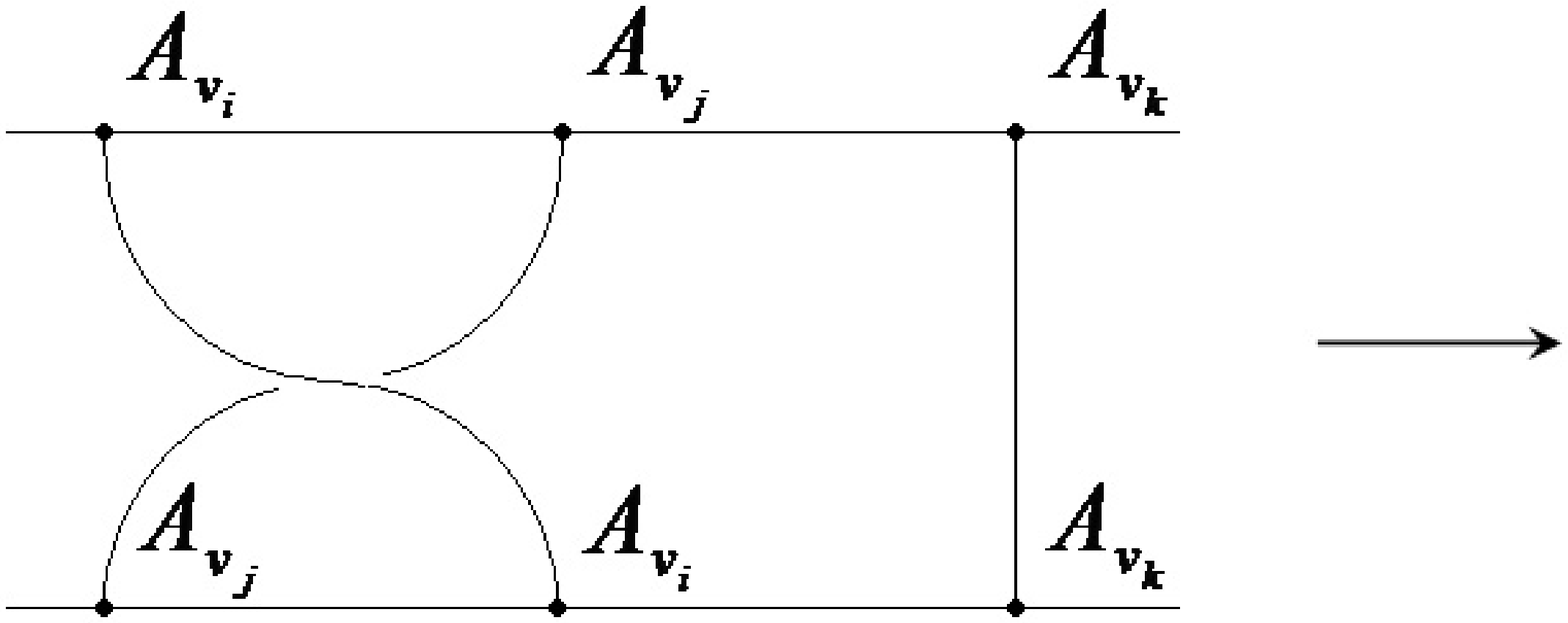}
\end{minipage}\hfill}
{\begin{minipage}{0.5\linewidth}$\Psi_{A_{v_{i}},A_{v_{j}}}\otimes
id_{A_{v_{k}}}$
\end{minipage}\hfill}
\end{figure}
%\vskip 1cm
\begin{figure}[!ht]
{\begin{minipage}{0.45\linewidth} \centering
\includegraphics[width=6cm,height=7.3cm]{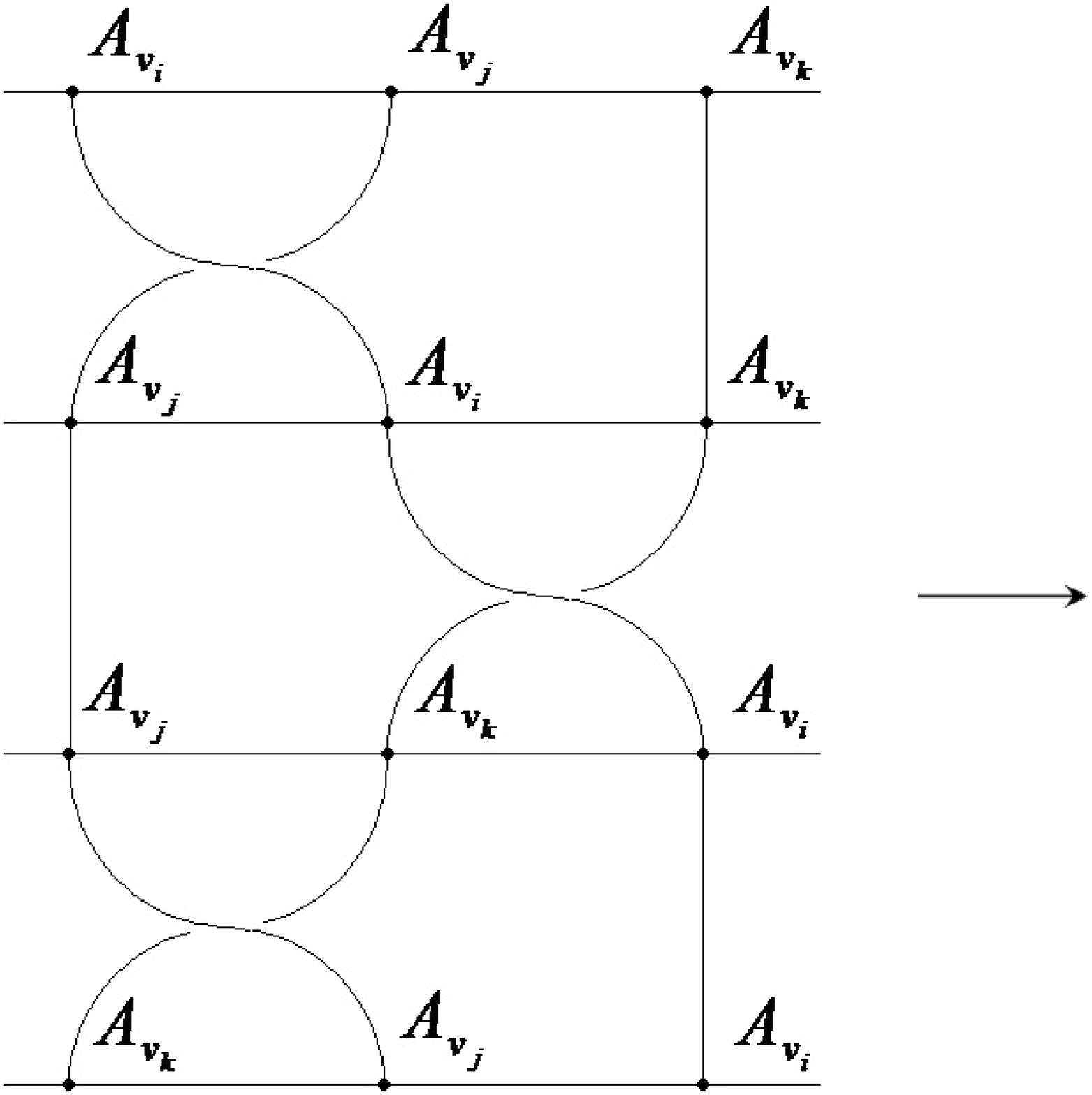}
\end{minipage}\hfill}
{\begin{minipage}{0.5\linewidth}$\left(\Psi_{A_{v_{j}},A_{v_{k}}}\otimes
id_{A_{v_{i}}}\right)\left(id_{A_{v_{j}}} \otimes
\Psi_{A_{v_{i}},A_{v_{k}}}\right)\\ \times
\left(\Psi_{A_{v_{i}},A_{v_{j}}}\otimes id_{A_{v_{k}}}\right)$
\end{minipage}\hfill}
\end{figure}

\vskip 1cm

\begin{figure}[!ht]
{\begin{minipage}{0.45\linewidth} \centering
\includegraphics[width=6cm,height=9.5cm]{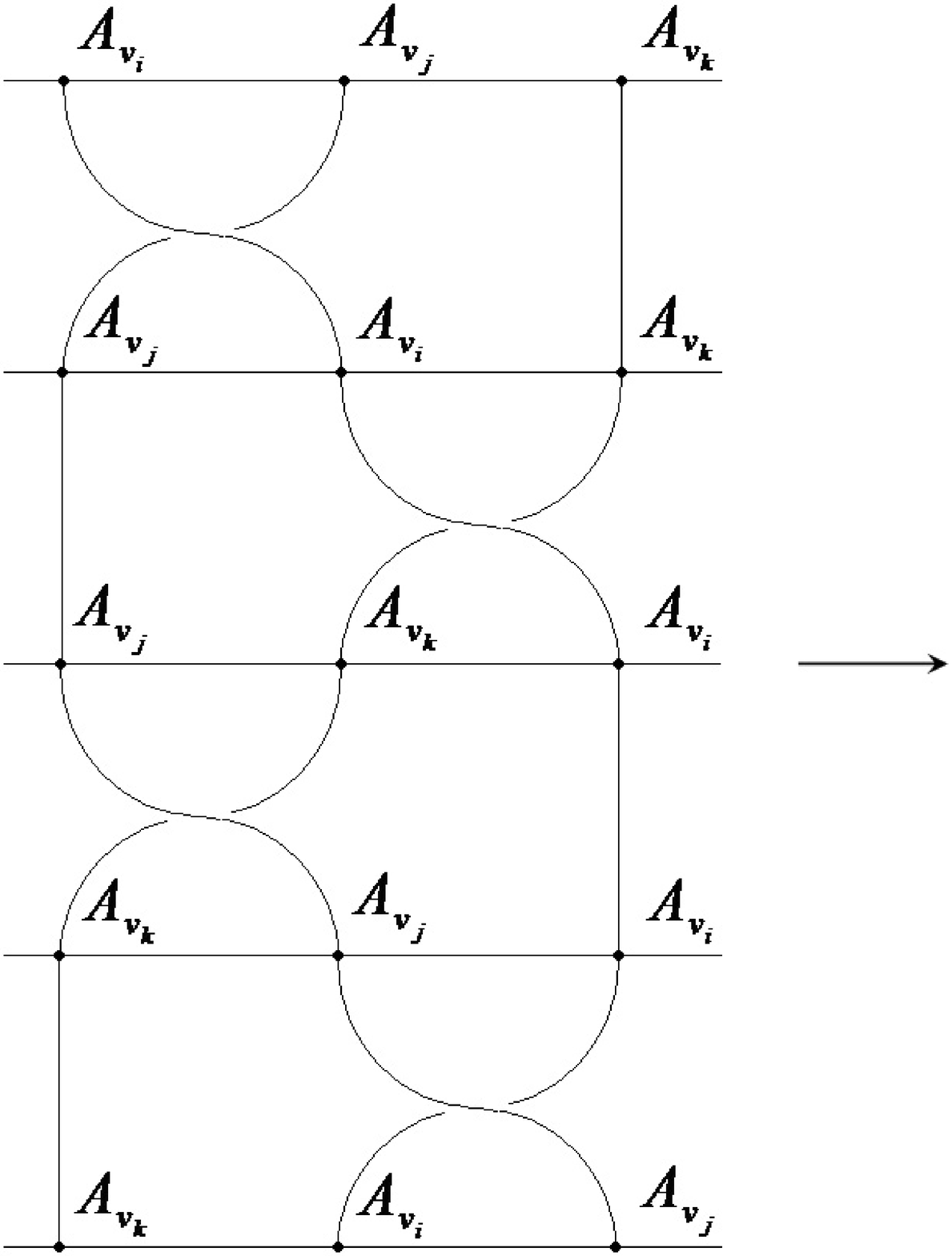}
\end{minipage}\hfill}
{\begin{minipage}{0.5\linewidth}$\left(id_{A_{v_{k}}}\otimes\Psi_{A_{v_{j}},A_{v_{i}}}\right)
\left(\Psi_{A_{v_{j}},A_{v_{k}}}\otimes id_{A_{v_{i}}}\right)
\times\left(id_{A_{v_{j}}} \otimes
\Psi_{A_{v_{i}},A_{v_{k}}}\right)\left(\Psi_{A_{v_{i}},A_{v_{j}}}\otimes
id_{A_{v_{k}}}\right)$
\end{minipage}\hfill}
\end{figure}

\newpage
%\vfill~\strut
\subsection{Right ternary maps}
By braiding from the right we construct the following maps
\vspace{0.4cm} \vskip 0.5cm
\begin{figure}[!ht]
{\begin{minipage}{0.45\linewidth} \centering
\includegraphics[width=6cm,height=2.7cm]{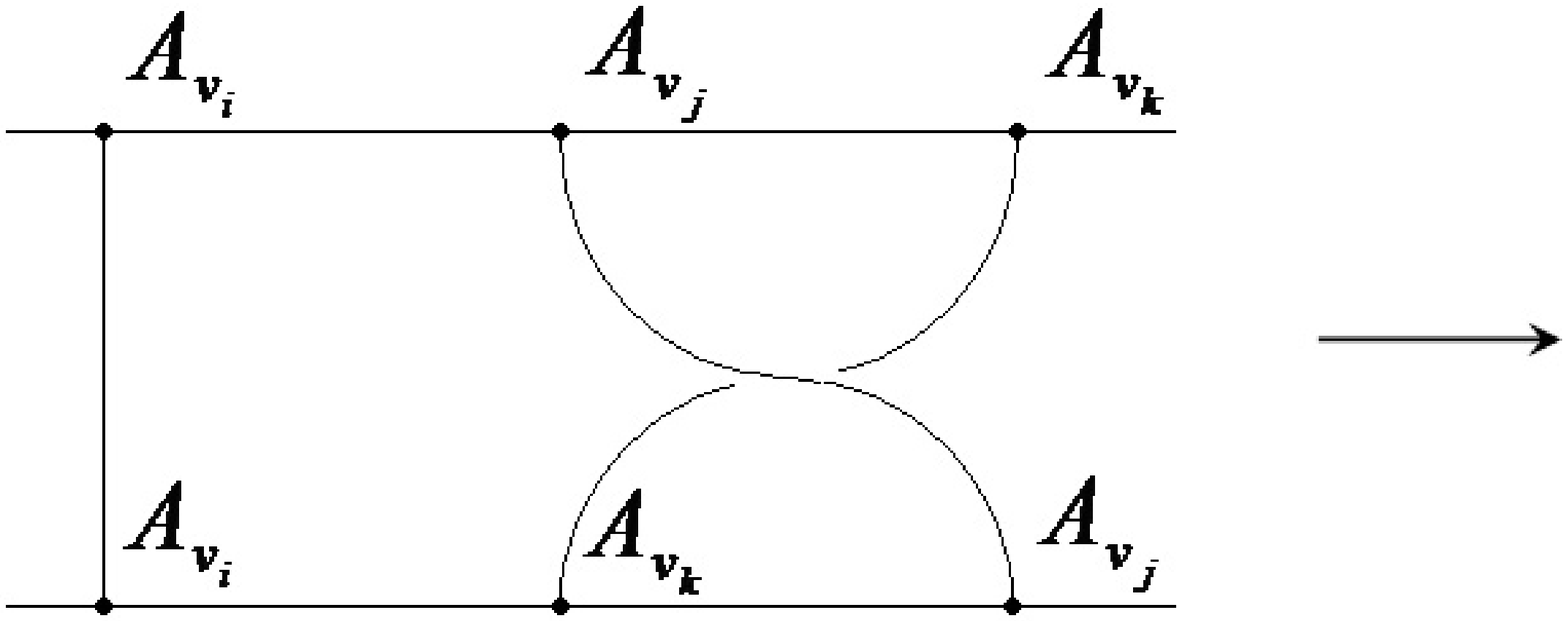}
\end{minipage}\hfill}
{\begin{minipage}{0.5\linewidth}$id_{A_{v_{i}}}\otimes
\Psi_{A_{v_{j}},A_{v_{k}}}$
\end{minipage}\hfill}
\end{figure}

\vskip 1cm
\begin{figure}[!h]
{\begin{minipage}{0.45\linewidth} \centering
\includegraphics[width=6cm,height=7.3cm]{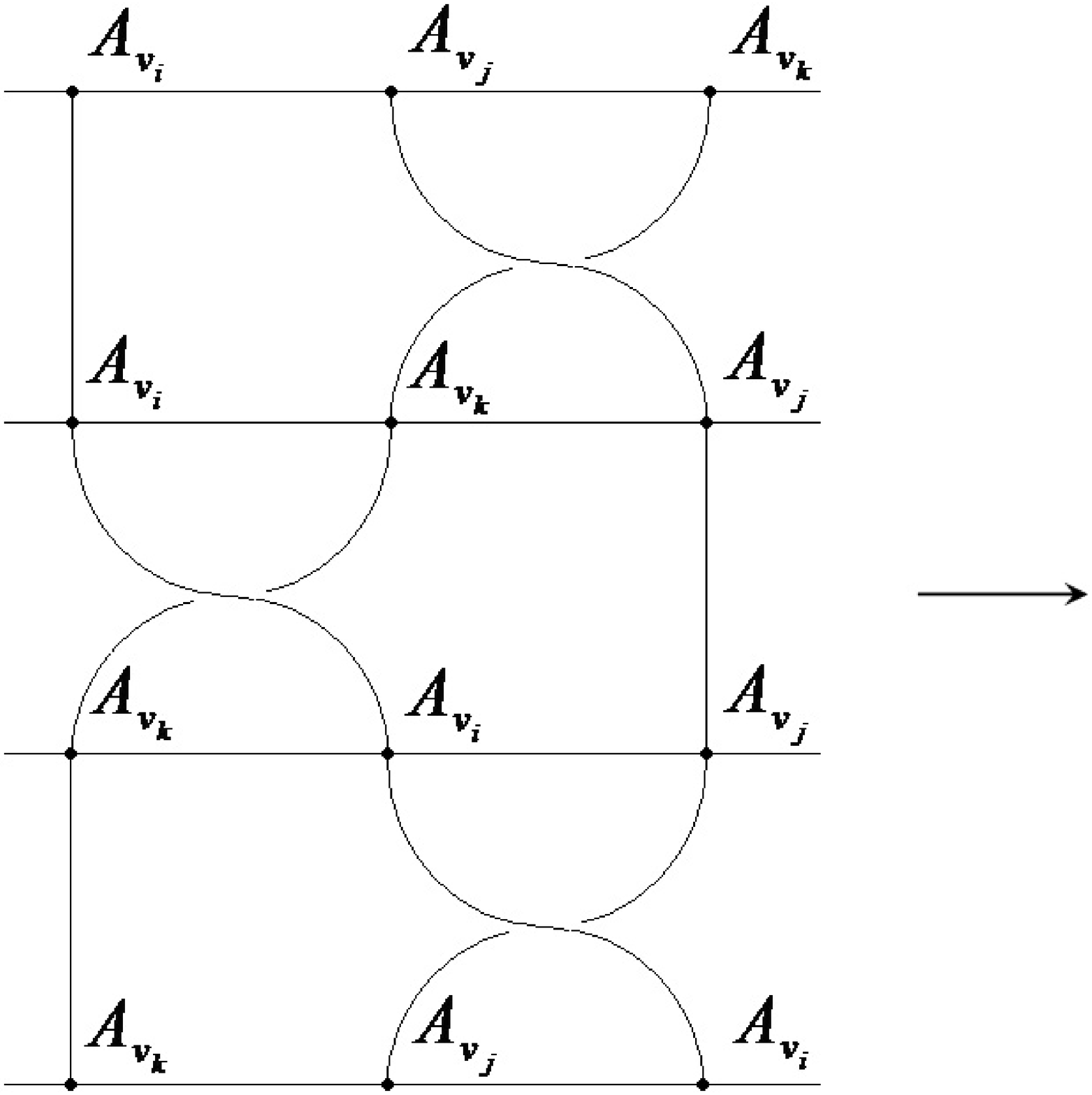}
\end{minipage}\hfill}
{\begin{minipage}{0.5\linewidth}$\left(id_{A_{v_{k}}} \otimes
\Psi_{A_{v_{i}},A_{v_{j}}}\right)
\left(\Psi_{A_{v_{i}},A_{v_{k}}}\otimes id_{A_{v_{j}}}\right)\\
\times\left(id_{A_{v_{i}}} \otimes
\Psi_{A_{v_{j}},A_{v_{k}}}\right)$
\end{minipage}\hfill}
\end{figure}

%\vfill~\strut
\vskip 1cm
\begin{figure}[!ht]
{\begin{minipage}{0.45\linewidth} \centering
\includegraphics[width=6cm,height=9.5cm]{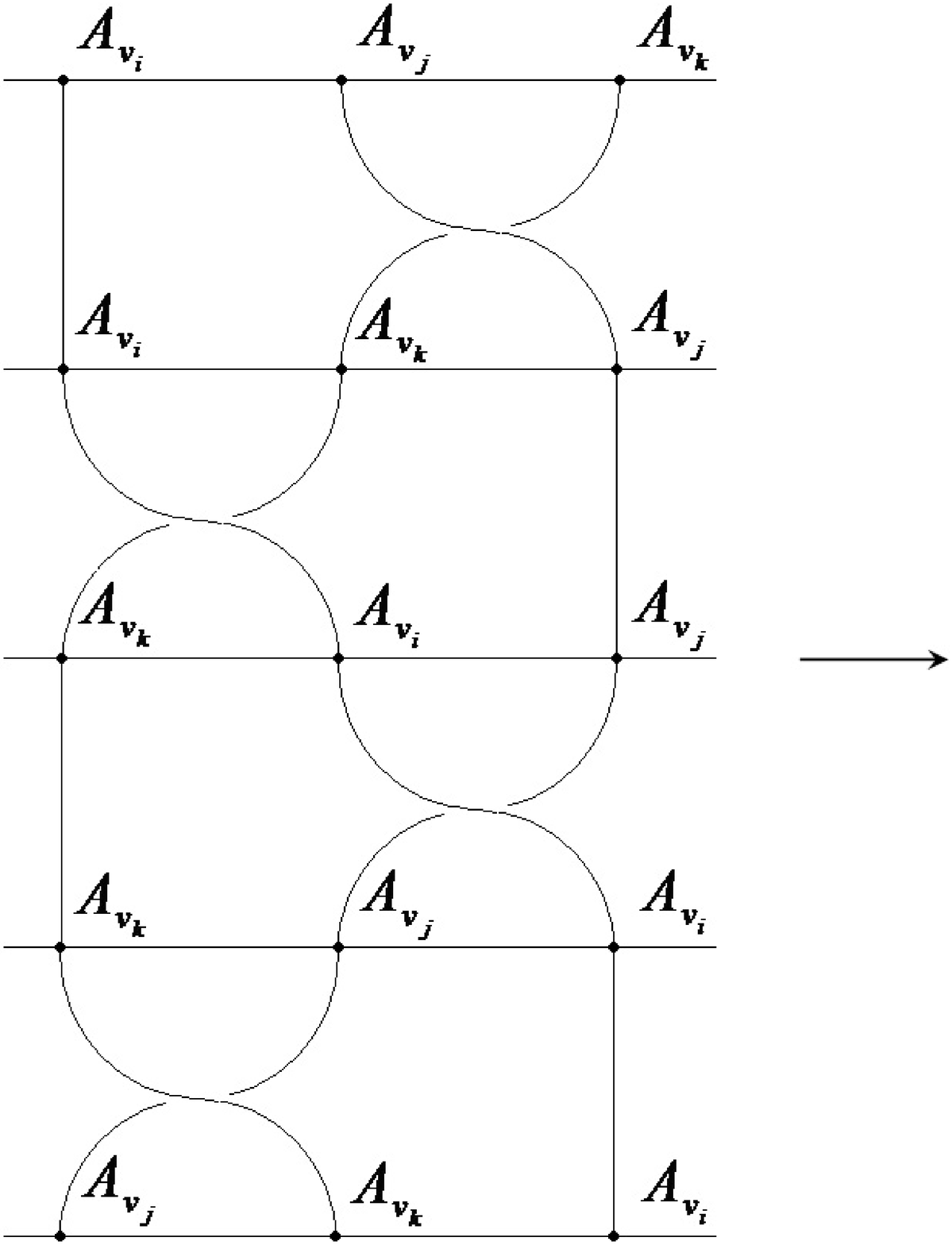}
\end{minipage}\hfill}
{\begin{minipage}{0.5\linewidth}$\left(\Psi_{A_{v_{k}},A_{v_{j}}}\otimes id_{A_{v_{i}}}\right)\left(id_{A_{v_{k}}} \otimes \Psi_{A_{v_{i}},A_{v_{j}}}\right)\\
\times\left(\Psi_{A_{v_{i}},A_{v_{k}}}\otimes
id_{A_{v_{j}}}\right)\left(id_{A_{v_{i}}} \otimes
\Psi_{A_{v_{j}},A_{v_{k}}}\right)$
\end{minipage}\hfill}
\end{figure}

\newpage
Let $\left<~,~ ,~\right>$ be a trilinear map defined on $A\times
A\times A$
\begin{eqnarray}
\left<~ ,~ ,~ \right>: A\times A\times A \longrightarrow A
\label{eq6}
\end{eqnarray}
Let $I_{1}$ be the two-sided ideal generated under the map
\begin{eqnarray}
& &id+\Psi_{A_{v_{i}},A_{v_{j}}}\otimes id_{A_{v_{k}}} \nonumber \\
& &-\left(id_{A_{v_{k}}}\otimes \Psi_{A_{v_{i}},A_{v_{j}}} \right)
\left(\Psi_{A_{v_{j}},A_{v_{k}}}\otimes id_{A_{v_{i}}} \right)
\left(id_{A_{v_{j}}}\otimes \Psi_{A_{v_{i}},A_{v_{k}}} \right)
\left(\Psi_{A_{v_{i}},A_{v_{j}}} \otimes id_{A_{v_{k}}} \right) \nonumber \\
& &-\left(\Psi_{A_{v_{j}},A_{v_{k}}}\otimes id_{A_{v_{i}}} \right)
\left(id_{A_{v_{j}}}\otimes \Psi_{A_{v_{i}},A_{v_{k}}} \right)
\left(\Psi_{A_{v_{i}},A_{v_{j}}} \otimes id_{A_{v_{k}}} \right)
-\left<~ ,~ ,~ \right> \label{eq7}
\end{eqnarray}
on $A_{v_{i}}\otimes A_{v_{j}}\otimes A_{v_{k}}$\\
Let $I_{2}$ be the two-sided ideal generated under the map
\begin{eqnarray}
& &id+id_{A_{v_{i}}}\otimes \Psi_{A_{v_{j}},A_{v_{k}}} \nonumber \\
& &-\left(\Psi_{A_{v_{k}},A_{v_{j}}}\otimes id_{A_{v_{i}}} \right)
\left(id_{A_{v_{k}}}\otimes \Psi_{A_{v_{i}},A_{v_{j}}} \right)
\left(\Psi_{A_{v_{i}},A_{v_{k}}} \otimes id_{A_{v_{j}}} \right)
\left(id_{A_{v_{i}}}\otimes \Psi_{A_{v_{j}},A_{v_{k}}} \right) \nonumber \\
& &-\left(id_{A_{v_{k}}}\otimes \Psi_{A_{v_{i}},A_{v_{j}}} \right)
\left(\Psi_{A_{v_{i}},A_{v_{k}}} \otimes id_{A_{v_{j}}} \right)
\left(id_{A_{v_{i}}}\otimes \Psi_{A_{v_{j}},A_{v_{k}}}\right)
-\left<~ ,~ ,~ \right> \label{eq8}
\end{eqnarray}
on $A_{v_{i}}\otimes A_{v_{j}}\otimes A_{v_{k}}$\\
where:\\
$id_{A_{v}}$ is the identity on the vector space $A_{v}$, that is
\begin{eqnarray}
id=id_{A_{v_{i}}}\otimes id_{A_{v_{j}}}\otimes id_{A_{v_{k}}}
\label{eq9}
\end{eqnarray}
\subsection{Envelopping algebras}
We define two enveloping algebras $\cup_{1}(A)=T(A)/I_{1}$ and
$\cup_{2}(A)=T(A)/I_{2}$. We denote by $\left<~,~ , ~\right>_{1}$
the trilinear map for $\cup_{1}(A)$ , and by $\left<~,~ ,
~\right>_{2}$ the trilinear map for $\cup_{2}(A)$. Left
$\left\{e_{i}\right\}$, $\left\{f_{j}\right\}$ and
$\left\{g_{k}\right\}$ be respectively the bases of the vector
spaces $A_{v_{i}}$, $A_{v_{j}}$ and $A_{v_{k}}$,
$\Psi=\left(\Psi_{ij}^{kl}\right)$ the matrix of
$\Psi_{A_{v_{i}},A_{v_{j}}}$ and denote the product in
$\cup_{1}(A)=T(A)/I_{1}$ and $\cup_{2}(A)=T(A)/I_{2}$ of elements
$x$ and $y$ in $A$ by $xy$, we have the following relations:\\

$\cup_{1}(A)=T(A)/I_{1}$ is an associative algebra; since the
composition of the canonical mapping of $T(A)$ onto $\cup_{1}(A)$
with the inclusion mapping of $A$ into $T(A)$ yields a one-to-one
mapping, we may identify $A$ with its image in $\cup_{1}(A)$\\
The trilinear map $\left<~,~ , ~\right>_{1}$ reads:
\begin{eqnarray}
\left<e_{i},f_{j},g_{k}\right>_{1}&=&e_{i}\otimes f_{j}\otimes g_{k}+\Psi_{ij}^{mn}f_{m}\otimes e_{n}\otimes g_{k} \nonumber \\
&&+\Psi_{ij}^{mn}\Psi_{nk}^{pq}\Psi_{mp}^{rs}\Psi_{sq}^{tu}g_{r}\otimes
e_{t}\otimes
f_{u}+\Psi_{ij}^{mn}\Psi_{nk}^{pq}\Psi_{mp}^{rs}g_{r}\otimes
f_{s}\otimes e_{q} \label{eq10}
\end{eqnarray}
Similarly for $\cup_{2}(A)$; we identity $A$ with its image in
$\cup_{2}(A)$.  The trilinear map  $\left<~,~ , ~\right>_{2}$ is
\begin{eqnarray}
\left<e_{i},f_{j},g_{k}\right>_{2}&=&e_{i}\otimes f_{j}\otimes g_{k}+\Psi_{jk}^{mn}e_{i}\otimes g_{m}\otimes f_{n} \nonumber \\
&&+\Psi_{jk}^{mn}\Psi_{im}^{pq}\Psi_{qn}^{rs}\Psi_{pr}^{tu}f_{t}\otimes
g_{u}\otimes
e_{s}+\Psi_{jk}^{mn}\Psi_{im}^{pq}\Psi_{qn}^{rs}g_{p}\otimes
f_{r}\otimes e_{s} \label{eq11}
\end{eqnarray}

\subsection{Definition of a para-algebra and Lie para-algebra :}
A para-algebra is a $\Gamma$-graded vector space
$A=\oplus_{i}A_{v_{i}}$; $v_{i}\in \Gamma$ (that is, if $a\in
A_{v_{i}}$, $b\in A_{v_{j}}$, $v_{i}$, $v_{j}\in \Gamma$, then
$ab\in A_{v_{i}+v_{j}}$); the braiding is taken to be
$\Psi_{A_{v_{i}},A_{v_{j}}}=(-1)^{\sigma(v_{i},v_{j})}$, where
$\sigma:\Gamma\times\Gamma\longrightarrow Z_{2}$,
$\left(Z_{2}=\left\{\overline{0},\overline{1}\right\}\right)$ is
defined in the paper$\cite {Zi}$,equipped with the following ternary
maps:
\begin{eqnarray}
\left<a,b,c\right>_{1}&=&a\otimes b\otimes c+(-1)^{\sigma(v_{i},v_{j})}b\otimes a\otimes c \nonumber \\
&&+(-1)^{\sigma(v_{i},v_{k})+\sigma(v_{j},v_{k})}c\otimes a\otimes
b+(-1)^{\sigma(v_{i},v_{j})+\sigma(v_{i},v_{k})+\sigma(v_{j},v_{k})}c\otimes
b\otimes a \label{eq12}
\end{eqnarray}
and
\begin{eqnarray}
\left<a,b,c\right>_{2}&=&a\otimes b\otimes c+(-1)^{\sigma(v_{j},v_{k})}a\otimes c\otimes b \nonumber \\
&&+(-1)^{\sigma(v_{i},v_{j})+\sigma(v_{i},v_{k})}b\otimes c\otimes
a+(-1)^{\sigma(v_{j},v_{k})+\sigma(v_{i},v_{j})+\sigma(v_{i},v_{k})}c\otimes
b\otimes a \label{eq13}
\end{eqnarray}
 A natural way of defining brackets $\left<~ ,~ ,
~\right>_{1}$ or $\left<~,~ , ~\right>_{2}$ in a para-algebra $A$;
is through the following equalities,
\begin{eqnarray}
\left<a,b,c\right>_{1}&=&(-1)^{\sigma(v_{i},v_{j})}\left<b,a,c\right>_{1} a\in A_{v_{i}}, b\in A_{v_{j}}, c\in A_{v_{k}} \label{eq14} \\
\left<a,b,c\right>_{2}&=&(-1)^{\sigma(v_{j},v_{k})}\left<a,c,b\right>_{2}
a\in A_{v_{i}}, b\in A_{v_{j}}, c\in A_{v_{k}} \label{eq15}
\end{eqnarray}
When we take $\left<~ ,~ , ~\right>_{1}$, we have the left
para-algebra, and when we take $\left<~,~ , ~\right>_{2}$, we have
the right para-algebra.\\
For an associative para-algebra $A$ the following identities hold:\\
\begin{eqnarray}
\left<a,b,cd\right>_{1}&=&\left<a,b,c\right>_{1}d+(-1)^{\sigma(v_{i}+v_{j},v_{k})}c\left<a,b,d\right>_{1} \label{eq16}\\
\left<a,b,cde\right>_{1}&=&\left<a,b,c\right>_{1}de+(-1)^{\sigma(v_{i}+v_{j},v_{k})}c\left<a,b,d\right>_{1}e \nonumber \\
& &
+(-1)^{\sigma(v_{i}+v_{j},v_{k}+v_{l})}cd\left<a,b,e\right>_{1} \label{eq17}\\
\left<ab,c,d\right>_{2}&=&a\left<b,c,d\right>_{2}+(-1)^{\sigma(v_{k}+v_{l},v_{j})}\left<a,c,d\right>_{2}b \label{eq18} \\
\left<abc,d,e\right>_{2}&=&ab\left<c,d,e\right>_{2}+(-1)^{\sigma(v_{k},v_{l}+v_{m})}a\left<b,d,e\right>_{2}c \nonumber \\
& & +(-1)^{\sigma(v_{j}+v_{k},v_{l}+v_{m})}\left<a,d,e\right>_{2}bc
\label{eq19}
\end{eqnarray}
where $a\in A_{v_{i}}$, $b\in A_{v_{j}}$, $c\in A_{v_{k}}$, $d\in
A_{v_{l}}$, $e\in A_{v_{m}}$\\
$A$ Lie para-algebra is a  para-algebra with a trilinear operation
$\left<~,~ , ~\right>_{1}$ or $\left<~,~ , ~\right>_{2}$ satisfying
the following relations:
\begin{eqnarray}
\left<a,b,c\right>_{1}=(-1)^{\sigma(v_{i},v_{j})}\left<b,a,c\right>_{1} a\in A_{v_{i}}, b\in A_{v_{j}}, c\in A_{v_{k}}  \label{eq20}\\
(-1)^{\sigma(v_{i},v_{k})}\left<a,b,c\right>_{1}+(-1)^{\sigma(v_{i},v_{j})}\left<b,c,a\right>_{1}+(-1)^{\sigma(v_{j},v_{k})}\left<c,a,b\right>_{1}=0
\label{eq21}\\
\left<a,b,\left<c,d,e\right>_{1}\right>_{1}=\left<\left<a,b,c\right>_{1},d,e\right>_{1}+(-1)^{\sigma(v_{i}+v_{j},v_{k})}\left<c,\left<a,b,d\right>_{1},e\right>_{1}
\nonumber \\
+(-1)^{\sigma(v_{i}+v_{j},v_{k})+\sigma(v_{i}+v_{j},v_{l})}\left<c,d,\left<a,b,e\right>_{1}\right>_{1}
\label{eq22}
\end{eqnarray}
where  $a\in A_{v_{i}}$, $b\in A_{v_{j}}$, $c\in A_{v_{k}}$, $d\in
A_{v_{l}}$, $e\in A_{v_{m}}$\\
The relations (\ref{eq20})- (\ref{eq22}) are the anticommutativity,
and the Jacobi identity for $\left<~,~ ,
~\right>_{1}$ respectively, while (\ref{eq21}) is a cyclic relation.\\
The relation (\ref{eq20}) - (\ref{eq22}) define
the left Lie para-algebra.\\
Similarly for $\left<~,~ , ~\right>_{2}$, we have
\begin{eqnarray}
\left<a,b,c\right>_{2}=(-1)^{\sigma(v_{j},v_{k})}\left<a,c,b\right>_{2} a\in A_{v_{i}}, b\in A_{v_{j}}, c\in A_{v_{k}}  \label{eq23}\\
(-1)^{\sigma(v_{i},v_{k})}\left<a,b,c\right>_{2}+(-1)^{\sigma(v_{i},v_{j})}\left<b,c,a\right>_{2}+(-1)^{\sigma(v_{j},v_{k})}\left<c,a,b\right>_{2}=0
\label{eq24}\\
\left<\left<a,b,c\right>_{2},d,e\right>_{2}=\left<a,b,\left<c,d,e\right>_{2}\right>_{2}+(-1)^{\sigma(v_{k},v_{l},v_{m})}\left<a,\left<b,d,e\right>_{2},c\right>_{2}
\nonumber \\
+(-1)^{\sigma(v_{j}+v_{k},v_{l}+v_{m})}\left<\left<a,d,e\right>_{2},b,c\right>_{2}
\label{eq25}
\end{eqnarray}
where  $a\in A_{v_{i}}$, $b\in A_{v_{j}}$, $c\in A_{v_{k}}$, $d\in
A_{v_{l}}$, $e\in A_{v_{m}}$\\
The relation (\ref{eq23}) - (\ref{eq25}) define a right Lie
para-algebra. The relations (\ref{eq23}), (\ref{eq25}) are
respectively the anticommutativity, and the Jacobi identity for
$\left<~,~ , ~\right>_{2}$, while (\ref{eq24}) is a cyclic relation.
The Jacobi identities (\ref{eq22}) and (\ref{eq25}) follow
respectively from (\ref{eq16}), (\ref{eq17}) and (\ref{eq18}),
(\ref{eq19}), respectively \\
We can define the trilinear maps $\left<~,~ , ~\right>_{1}$ and
$\left<~,~ , ~\right>_{2}$ by introducing respectively the left
multiplication operator and the right multiplication operator which
we denote by $L_{a,b}$ and $R_{y,z}$ such that:
\begin{eqnarray}
L_{a,b}&:& A\longrightarrow A
\nonumber \\
L_{a,b}&:& c\longrightarrow
L_{a,b}c=\left<a,b,c\right>_{1},a,b~~~~\mbox{and}~~~~c\in A
\label{eq26}
\end{eqnarray}
and,
\begin{eqnarray}
R_{y,z}&:& A\longrightarrow A
\nonumber \\
R_{y,z}&:& x\longrightarrow
R_{y,z}(x)=xR_{y,z}=\left<x,y,z\right>_{2},x,y~~~~\mbox{and}~~~~z\in
A \label{eq27}
\end{eqnarray}
The operator $L_{a,b}$ behave just like  first order differential
operator in that it obeys a product rule
\begin{eqnarray}
L_{a,b}\left<c,d,e\right>_{1}&=&\left<L_{a,b}c,d,e\right>_{1}+(-1)^{\sigma(v_{i}+v_{j},v_{k})}\left<c,L_{a,b}d,e\right>_{1}
\nonumber \\
&&+(-1)^{\sigma(v_{i}+v_{j},v_{k})+\sigma(v_{i}+v_{j},v_{l})}\left<c,d,L_{a,b}e\right>_{1}
\label{eq28}
\end{eqnarray}
$a\in A_{v_{i}}$, $b\in A_{v_{j}}$, $c\in A_{v_{k}}$, $d\in
A_{v_{l}}$, $e\in A_{v_{m}}$\\
This product rule is simply another way of writing the
Jacobi identity (\ref{eq22}) for $\left<~,~ , ~\right>_{1}$.\\
The operator $R_{y,z}$ has the property
\begin{eqnarray}
\left<a,b,c\right>_{2}R_{d,e}&=&\left<a,b,cR_{d,e}\right>_{2}+(-1)^{\sigma(v_{k},v_{l}+v_{m})}\left<a,bR_{d,e},c\right>_{2}
\nonumber \\
&&+(-1)^{\sigma(v_{j}+v_{k},v_{l}+v_{m})}\left<aR_{d,e},b,c\right>_{2}
\label{eq29}
\end{eqnarray}
for $a\in A_{v_{i}}$, $b\in A_{v_{j}}$, $c\in A_{v_{k}}$, $d\in
A_{v_{l}}$, $e\in A_{v_{m}}$. This is also another way of writing
the Jacobi identity (\ref{eq25}) for $\left<~,~ , ~\right>_{2}$.

\section{Application to the formalism of parastatistics.}

Since $A$ is a vector space, we define a quadratic  form $Q$ on $A$,
which takes its values in the field $K$, such that:
\begin{eqnarray}
Q(a,b)=-(-1)^{\sigma(v_{i},v_{j})}Q(b,a),a\in A_{v_{i}}, b\in
A_{v_{j}}.  \label{eq30}
\end{eqnarray}
Let $J_{1}$ be the two sided ideal in $T(A)$ generated by the
elements:
\begin{eqnarray}
&&a\otimes b\otimes c+(-1)^{\sigma(v_{i},v_{j})}b\otimes a\otimes
c-(-1)^{\sigma(v_{i},v_{k})+\sigma(v_{j},v_{k})}c\otimes a\otimes b
\nonumber \\
&&-(-1)^{\sigma(v_{i},v_{j})+\sigma(v_{i},v_{k})+\sigma(v_{j},v_{k})}c\otimes
b\otimes a \nonumber \\
&&-\left\{aQ(b,c)+(-1)^{\sigma(v_{i},v_{j})}bQ(a,c)-(-1)^{\sigma(v_{i},v_{k})+\sigma(v_{j},v_{k})}Q(c,a)b\right. \nonumber \\
&&\left.-(-1)^{\sigma(v_{i},v_{j})+\sigma(v_{i},v_{k})+\sigma(v_{j},v_{k})}Q(c,b)a\right\}.
\label{eq31}
\end{eqnarray}
where $a\in A_{v_{i}}$, $b\in A_{v_{j}}$, $c\in A_{v_{k}}$.Let
$J_{2}$ be the two sided ideal in $T(A)$ generated by the elements:
\begin{eqnarray}
&&a\otimes b\otimes c+(-1)^{\sigma(v_{j},v_{k})}a\otimes c\otimes
b-(-1)^{\sigma(v_{i},v_{j})+\sigma(v_{i},v_{k})}b\otimes c\otimes a
\nonumber \\
&&-(-1)^{\sigma(v_{j},v_{k})+\sigma(v_{i},v_{j})+\sigma(v_{i},v_{k})}c\otimes
b\otimes a \nonumber \\
&&-\left\{Q(a,b)c+(-1)^{\sigma(v_{j},v_{k})}Q(a,c)b-(-1)^{\sigma(v_{i},v_{j})+\sigma(v_{i},v_{k})}bQ(c,a)\right. \nonumber \\
&&\left.-(-1)^{\sigma(v_{j},v_{k})+\sigma(v_{i},v_{j})+\sigma(v_{i},v_{k})}cQ(b,a)\right\}.
\label{eq32}
\end{eqnarray}
where $a\in A_{v_{i}}$, $b\in A_{v_{j}}$, $c\in A_{v_{k}}$.\\
We define the envelopping  algebra $V_{1}(A)=T(A)/J_{1}$ and
$V_{2}(A)=T(A)/J_{2}$; $V_{1}(A)$ and $V_{2}(A)$ are associative
algebras, we denote the product of two elements $x$ and $y$ in
$V_{1}(A)$ or in $V_{2}(A)$ by
$xy$.\\
Since $A$ is embeded in $V_{1}(A)$ then in $A$ we have
\begin{eqnarray}
&&abc+(-1)^{\sigma(v_{i},v_{j})}bac-(-1)^{\sigma(v_{i},v_{k})+\sigma(v_{j},v_{k})}cab-(-1)^{\sigma(v_{i},v_{j})+\sigma(v_{i},v_{k})+\sigma(v_{j},v_{k})}cba= \nonumber \\
&&aQ(b,c)+(-1)^{\sigma(v_{i},v_{j})}bQ(a,c)-(-1)^{\sigma(v_{i},v_{k})+\sigma(v_{j},v_{k})}Q(c,a)b \nonumber \\
&&-(-1)^{\sigma(v_{i},v_{j})+\sigma(v_{i},v_{k})+\sigma(v_{j},v_{k})}Q(c,b)a
\label{eq33}
\end{eqnarray}
Similarly in $A$ we have
\begin{eqnarray}
&&abc+(-1)^{\sigma(v_{j},v_{k})}bac-(-1)^{\sigma(v_{i},v_{j})+\sigma(v_{i},v_{k})}cab-(-1)^{\sigma(v_{j},v_{k})+\sigma(v_{i},v_{j})+\sigma(v_{i},v_{k})}cba= \nonumber \\
&&Q(a,b)c+(-1)^{\sigma(v_{j},v_{k})}Q(a,c)b-(-1)^{\sigma(v_{i},v_{j})+\sigma(v_{i},v_{k})}bQ(c,a)b \nonumber \\
&&-(-1)^{\sigma(v_{j},v_{k})+\sigma(v_{i},v_{j})+\sigma(v_{i},v_{k})}cQ(b,a)
\label{eq34}
\end{eqnarray}
If we compare (\ref{eq12}) and (\ref{eq33}), then the trilinear map
$\left<~,~ , ~\right>_{1}$ reads in this case
\begin{eqnarray}
\left<a,b,c\right>_{1}&=&aQ(b,c)+(-1)^{\sigma(v_{i},v_{j})}bQ(q,c)-(-1)^{\sigma(v_{i},v_{k})+\sigma(v_{j},v_{k})}Q(c,a)b \nonumber \\
&&-(-1)^{\sigma(v_{i},v_{j})+\sigma(v_{i},v_{k})+\sigma(v_{j},v_{k})}Q(c,b)a
\label{eq35}
\end{eqnarray}
Comparing (\ref{eq13}) and (\ref{eq34}) shows in this case the
trilinear map $\left<~,~ , ~\right>_{2}$ reads
\begin{eqnarray}
\left<a,b,c\right>_{2}&=&Q(a,b)c+(-1)^{\sigma(v_{j},v_{k})}Q(a,c)b-(-1)^{\sigma(v_{i},v_{j})+\sigma(v_{i},v_{k})}bQ(c,a) \nonumber \\
&&-(-1)^{\sigma(v_{j},v_{k})+\sigma(v_{i},v_{j})+\sigma(v_{i},v_{k})}cQ(b,a)
\label{eq36}
\end{eqnarray}

\subsection{Exemple}

when $\Gamma$ is a two-dimensional vector space on $Z_{2}$,
$\Gamma=\left\{\left(\begin{array}{c}0
\\
0\end{array}\right), \left(\begin{array}{c}1
\\
0\end{array}\right), \left(\begin{array}{c}0
\\
1\end{array}\right), \left(\begin{array}{c}1
\\
1\end{array}\right)\right\}$, the algebra $A$ is $
A=A_{\left(\begin{array}{c}0
\\
0\end{array}\right)}\oplus A_{\left(\begin{array}{c}1
\\
0\end{array}\right)}\oplus A_{\left(\begin{array}{c}0
\\
1\end{array}\right)}\oplus A_{\left(\begin{array}{c}1
\\
1\end{array}\right)}\\
$
 In the paper$^{5}$ we have shown that the
color superalgebra $C(2,s)$, the Lie superalgebra $C(1,s)$, the
color algebra $C(2,a)$ and the Lie algebra are caracterized,
respectively, by the following four  equivalence classes,
\begin{eqnarray}
\left\{\left(\begin{array}{cc} 1 & 0
\\
0 & 1 \end{array}\right), \left(\begin{array}{cc} 1 & 1
\\
1 & 0 \end{array}\right), \left(\begin{array}{cc} 0 & 1
\\
1 & 1 \end{array}\right)\right\}. \label{eq37}
\end{eqnarray}
\begin{eqnarray}
\left\{\left(\begin{array}{cc} 1 & 0
\\
0 & 0 \end{array}\right), \left(\begin{array}{cc} 0 & 0
\\
0 & 1 \end{array}\right), \left(\begin{array}{cc} 1 & 1
\\
1 & 1 \end{array}\right)\right\}. \label{eq38}
\end{eqnarray}
\begin{eqnarray}
\left\{\left(\begin{array}{cc} 0 & 1
\\
1 & 0 \end{array}\right)\right\}. \label{eq39}
\end{eqnarray}
\begin{eqnarray}
\left\{\left(\begin{array}{cc} 0 & 0
\\
0 & 0 \end{array}\right)\right\}. \label{eq40}
\end{eqnarray}
While $\sigma$ is represented by $M\in Sbil(Z_{2})$$^{5}$, such that:\\
$\sigma(v_{i},v_{j})=v_{i}^{t}Mv_{j}$, $v_{i}$, $v_{j}\in \Gamma$
and $v_{i}^{t}$ is the transpose of $v_{i}$.\\
Let $M=\left(\begin{array}{cc} \alpha_{11} & \alpha_{12}
\\
\alpha_{21} & \alpha_{22} \end{array}\right)$,
$v_{i}=\left(\begin{array}{c}p_{1}
\\
p_{2}\end{array}\right)$ and $v_{j}=\left(\begin{array}{c}q_{1}
\\
q_{2}\end{array}\right)$ where $\alpha_{11}$, $\alpha_{12}$,
$\alpha_{21}$, $\alpha_{22}$, $p_{1}$, $q_{1}$, $p_{2}$ and
$q_{2}\in Z_{2}$
\begin{eqnarray}
\sigma(v_{i},v_{j})=v_{i}^{t}Mv_{j}=p_{1}\alpha_{11}q_{1}+p_{1}\alpha_{12}q_{2}+p_{2}\alpha_{21}q_{1}+p_{2}\alpha_{22}q_{2}
\label{eq41}
\end{eqnarray}
Let:\\
$\left\langle E \right\rangle$ ($E$ is the identity) be the
generator of $A_{\left(\begin{array}{c}0
\\
0\end{array}\right)}$;\\
$\left\langle a_{i},a_{j}^{+},i,j=1,2,...,n\right\rangle$ be the
generator of $A_{\left(\begin{array}{c}1
\\
0\end{array}\right)}$;\\
$\left\langle b_{i},b_{j}^{+},i,j=1,2,...,m\right\rangle$ be the
generator of $A_{\left(\begin{array}{c}0
\\
1\end{array}\right)}$;\\
$\left\langle c_{i},c_{j}^{+},i,j=1,2,...,p\right\rangle$ be the
generator of $A_{\left(\begin{array}{c}1
\\
1\end{array}\right)}$.\\
and $M=\left(\begin{array}{cc} \alpha_{11} & \alpha_{12}
\\
\alpha_{21} & \alpha_{22}
\end{array}\right)=\left(\begin{array}{cc} 1 &
0
\\
0 & 1 \end{array}\right)$.\\
$A_{\left(\begin{array}{c}0
\\
0\end{array}\right)}$ and $A_{\left(\begin{array}{c}1
\\
1\end{array}\right)}$ are even spaces (symmetric spaces or
equivalently bosonic spaces) with respect to $\sigma$ since
$\sigma\left(\left(\begin{array}{c}0
\\
0\end{array}\right), \left(\begin{array}{c}0
\\
0\end{array}\right)\right) = 0$ and
$\sigma\left(\left(\begin{array}{c}1
\\
1\end{array}\right), \left(\begin{array}{c}1
\\
1\end{array}\right)\right)=\left(\begin{array}{cc} 1 &
1\end{array}\right)\left(\begin{array}{cc} 1 & 0
\\
0 & 1 \end{array}\right)\left(\begin{array}{c}1
\\
1\end{array}\right)=2=0$, with implies that if $a$ and $b\in
A_{\left(\begin{array}{c}0
\\
0\end{array}\right)}$ or $a$ and $b\in A_{\left(\begin{array}{c}1
\\
1\end{array}\right)}$ then $Q(a,b)=-Q(b,a)$.\\
$A_{\left(\begin{array}{c}1
\\
0\end{array}\right)}$ and $A_{\left(\begin{array}{c}0
\\
1\end{array}\right)}$ are odd spaces (antisymmetric spaces or
equivalently fermionic spaces) with respect to $\sigma$ since
$\sigma\left(\left(\begin{array}{c}1
\\
0\end{array}\right), \left(\begin{array}{c}1
\\
0\end{array}\right)\right)=\left(\begin{array}{cc} 1 &
0\end{array}\right)\left(\begin{array}{cc} 1 & 0
\\
0 & 1 \end{array}\right)\left(\begin{array}{c}1
\\
0\end{array}\right)=1$ and $\sigma\left(\left(\begin{array}{c}0
\\
1\end{array}\right), \left(\begin{array}{c}0
\\
1\end{array}\right)\right)=\left(\begin{array}{cc} 0 &
1\end{array}\right)\left(\begin{array}{cc} 1 & 0
\\
0 & 1 \end{array}\right)\left(\begin{array}{c}0
\\
1\end{array}\right)=1$, then $Q(a,b)=Q(b,a)$ for $a$ and $b\in
A_{\left(\begin{array}{c}1
\\
0\end{array}\right)}$ or $a$ and $b\in A_{\left(\begin{array}{c}0
\\
1\end{array}\right)}$.\\
First we consider the case where $A=A_{\left(\begin{array}{c}0
\\
0\end{array}\right)}\oplus A_{\left(\begin{array}{c}1
\\
0\end{array}\right)}$, $A_{\left(\begin{array}{c}0
\\
1\end{array}\right)}=\emptyset$ and $A_{\left(\begin{array}{c}1
\\
1\end{array}\right)}=\emptyset$\\

$a$, $b$ and $c$ can be taken as annihilation or creation operators,
$a_{i}$ and $a_{j}^{+}$ $i,j=1,2,...,n$. The quadratic form is
symmetric $Q(a,b)=Q(b,a)$; we denote the usual commutator by
$\left[a,b\right]=ab-ba$ and the usual anticommutator by $\left\{
a,b\right\}=ab+ba$.
\begin{eqnarray}
\left\langle
a,b,c\right\rangle_{1}&=&aQ(b,c)-bQ(a,c)-Q(c,a)b+Q(c,b)a=2Q(b,c)a-2Q(a,c)b
\nonumber \\
&=&abc-bac-cab+cba \nonumber \\
&=&\left[\left[a,b\right],c\right]  \label{eq42}
\end{eqnarray}
if we set, $a=a_{i}$, $b=a_{j}^{+}$ and $c=a_{k}$, and choose the
bilinear form $Q(a,b)$ to satisfy the following relations:
$Q(a_{i},a_{j}^{+})=Q(a_{j}^{+},a_{i})=\delta_{ij}$,
$Q(a_{i},a_{j}^{+})=Q(a_{j}^{+},a_{i})=\delta_{ij}$,
$Q(a_{i}^{+},a_{j}^{+})=Q(a_{j}^{+},a_{i}^{+})=0$,
$Q(a_{i},a_{j})=Q(a_{j},a_{i})=0$,then
\begin{eqnarray}
\left\langle
a_{i},a_{j}^{+},a_{k}\right\rangle_{1}&=&a_{i}Q(a_{j}^{+},a_{k})-Q(a_{i},a_{k})a_{j}^{+}-a_{j}^{+}Q(a_{k},a_{i})+Q(a_{k},a_{j}^{+})a_{i} \nonumber \\
&=&a_{i}\delta_{jk}+\delta_{jk}a_{i}=2\delta_{jk}a_{i}=\left[\left[a_{i},a_{j}^{+}\right],a_{k}\right]
\label{eq43}
\end{eqnarray}
when $a=a_{i}^{+}$, $b=a_{j}^{+}$ and $c=a_{k}$, it readily follows
that
\begin{eqnarray}
\left\langle
a_{i}^{+},a_{j}^{+},a_{k}\right\rangle_{1}=2\delta_{jk}a_{i}^{+}-2\delta_{ik}a_{j}^{+}=\left[\left[a_{i}^{+},a_{j}^{+}\right],a_{k}\right]
\label{eq44}
\end{eqnarray}
Eq. (\ref{eq44}) and Eq. (\ref{eq45}) are two of the relations that
caracterize the parafermion statistic. The others relations
are obtained in a similar way\\

\subsection{Remark}

One can take the trilinear application $\left\langle ~,~ ,~
\right\rangle_{2}$
\begin{eqnarray}
\left\langle
a_{i},a_{j}^{+},a_{k}\right\rangle_{2}=+\delta_{ij}a_{k}+a_{k}\delta_{ij}=2\delta_{ij}a_{k}=\left[a_{i},\left[a_{j}^{+},a_{k}\right]\right]
\label{eq45}
\end{eqnarray}
For $a=a_{i}^{+}$, $b=a_{j}^{+}$ and $c=a_{k}$
\begin{eqnarray}
\left\langle
a_{i}^{+},a_{j}^{+},a_{k}\right\rangle_{2}=-2\delta_{ik}a_{j}^{+}=\left[a_{i},\left[a_{j}^{+},a_{k}\right]\right]
\label{eq46}
\end{eqnarray}
the others relations follow directly.\\
Therefore, the present example reproduces the parafermion
statistics.\\
Now let $A=A_{\left(\begin{array}{c}0
\\
0\end{array}\right)}\oplus A_{\left(\begin{array}{c}1
\\
1\end{array}\right)}$, $\sigma\left(\left(\begin{array}{c}1
\\
1\end{array}\right), \left(\begin{array}{c}1
\\
1\end{array}\right)\right)=0$.\\
$a$, $b$ and $c\in A_{\left(\begin{array}{c}1
\\
1\end{array}\right)}$ may be identified with the annihilation and
creation operators $c_{i}$ and $c_{j}^{+}$ $i,j=1,2,...,n$, the
quadratic form satisfies $Q(a,b)=-Q(b,a)$ such that:
\begin{eqnarray}
Q(c_{j}^{+},c_{i})=-Q(c_{i},c_{j}^{+})=\delta_{ij}~~~~\mbox{and}~~~~Q(c_{i},c_{j})=-Q(c_{i}^{+},c_{j}^{+})=0
\label{eq47}
\end{eqnarray}
If we set $a=c_{i}$, $b=c_{j}^{+}$ and $c=c_{k}$ then:
\begin{eqnarray}
\left\langle
c_{i},c_{j}^{+},c_{k}\right\rangle_{1}&=&c_{i}Q(c_{j}^{+},c_{k})+c_{j}^{+}Q(c_{i},c_{k})-Q(c_{k},c_{i})c_{j}^{+}-Q(c_{k},c_{j}^{+})c_{i} \nonumber \\
&=&c_{i}\delta_{jk}+\delta_{jk}c_{i}=2\delta_{jk}c_{i} \nonumber \\
&=&c_{i}c_{j}^{+}c_{k}+c_{j}^{+}c_{i}c_{k}-c_{k}c_{i}c_{j}^{+}-c_{k}c_{j}^{+}c_{i}
\label{eq48}
\end{eqnarray}
In term of the usual commutator and anticommutator the expression (48) is  $\left[\left\{c_{i},c_{j}^{+}\right\},c_{k}\right]$\\
If we set $a=c_{i}$, $b=c_{j}^{+}$ and $c=c_{k}^{+}$
\begin{eqnarray}
\left\langle
c_{i},c_{j}^{+},c_{k}^{+}\right\rangle_{1}&=&c_{i}Q(c_{j}^{+},c_{k}^{+})+c_{j}^{+}Q(c_{i},c_{k}^{+})-Q(c_{k}^{+},c_{i})c_{j}^{+}-Q(c_{k}^{+},c_{j}^{+})c_{i} \nonumber \\
&=&-c_{j}^{+}\delta_{ik}-\delta_{ik}c_{j}^{+}=-2\delta_{ik}c_{j}^{+} \nonumber \\
&=&c_{i}c_{j}^{+}c_{k}^{+}+c_{j}^{+}c_{i}c_{k}^{+}-c_{k}^{+}c_{i}c_{j}^{+}-c_{k}^{+}c_{j}^{+}c_{i}  \nonumber \\
&=&\left[\left\{ c_{i},c_{j}^{+}\right\},c_{k}^{+}\right]
\label{eq49}
\end{eqnarray}
If we set $a=c_{i}$, $b=c_{j}^{+}$ and $c=c_{k}^{+}$ a direct
calculation gives,
\begin{eqnarray}
\left\langle
c_{i},c_{j},c_{k}^{+}\right\rangle_{1}&=&-2\delta_{jk}c_{i}-2\delta_{ik}c_{j} \nonumber \\
&=&\left[\left\{ c_{i},c_{j}\right\},c_{k}^{+}\right] \label{eq50}
\end{eqnarray}
If we set $a=c_{i}^{+}$, $b=c_{j}^{+}$ and $c=c_{k}$ then,
\begin{eqnarray}
\left\langle
c_{i}^{+},c_{j}^{+},c_{k}\right\rangle_{1}&=&2\delta_{jk}c_{i}^{+}+2\delta_{ik}c_{j}^{+}=\left[\left\{
c_{i}^{+},c_{j}^{+}\right\},c_{k}\right]
\label{eq51}\\
\left\langle c_{i},c_{j},c_{k}\right\rangle_{1}&=&0=\left[\left\{
c_{i},c_{j}\right\},c_{k}\right]\label{eq52}
\end{eqnarray}
The case $a=c_{i}^{+}$, $b=c_{j}^{+}$ and $c=c_{k}^{+}$
\begin{eqnarray}
\left\langle
c_{i}^{+},c_{j}^{+},c_{k}^{+}\right\rangle_{1}=0=\left[\left\{
c_{i}^{+},c_{j}^{+}\right\},c_{k}^{+}\right]\label{eq53}
\end{eqnarray}
therefore reproduces the paraboson statistics.\\
One can use the map $\left\langle ~,~ ,~ \right\rangle_{2}$, we
found that:
\begin{eqnarray}
\left\langle
c_{i},c_{j}^{+},c_{k}\right\rangle_{2}&=&Q(c_{i},c_{j}^{+})c_{k}+Q(c_{i},c_{k})c_{j}^{+}-c_{j}^{+}Q(c_{k},c_{i})-c_{k}Q(c_{j}^{+},c_{i}) \nonumber \\
&=&-\delta_{ij}c_{k}-c_{k}\delta_{ij}=-2\delta_{ij}c_{k} \nonumber \\
&=&c_{i}c_{j}^{+}c_{k}+c_{i}c_{k}c_{j}^{+}-c_{j}^{+}c_{k}c_{i}-c_{k}c_{j}^{+}c_{i}  \nonumber \\
&=&\left[c_{i}\left\{c_{j}^{+},c_{k}\right\}\right] \label{eq54}
\end{eqnarray}
\begin{eqnarray}
\left\langle
c_{i},c_{j}^{+},c_{k}^{+}\right\rangle_{2}&=&Q(c_{i},c_{j}^{+})c_{k}^{+}+Q(c_{i},c_{k}^{+})c_{j}^{+}-c_{j}^{+}Q(c_{k}^{+},c_{i})-c_{k}^{+}Q(c_{j}^{+},c_{i}) \nonumber \\
&=&-\delta_{ij}c_{k}^{+}-\delta_{ik}c_{j}^{+}-c_{j}^{+}\delta_{ik}-c_{k}^{+}\delta_{ij}=-2\delta_{ij}c_{k}^{+}-2\delta_{ik}c_{j}^{+} \nonumber \\
&=&c_{i}c_{j}^{+}c_{k}^{+}+c_{j}^{+}c_{i}c_{k}^{+}-c_{k}^{+}c_{i}c_{j}^{+}-c_{k}^{+}c_{j}^{+}c_{i}  \nonumber \\
&=&\left[\left\{ c_{i},c_{j}^{+}\right\},c_{k}^{+}\right]\nonumber \\
&=&\left[c_{i}\left\{c_{j}^{+},c_{k}^{+}\right\}\right] \label{eq55}
\end{eqnarray}
A similar calcultion gives;
\begin{eqnarray}
\left\langle
c_{i},c_{j},c_{k}^{+}\right\rangle_{2}&=&-2\delta_{ik}c_{j}=\left[c_{i}\left\{c_{j},c_{k}^{+}\right\}\right]
\label{eq56}\\
\left\langle
c_{i}^{+},c_{j}^{+},c_{k}\right\rangle_{2}&=&2\delta_{ik}c_{j}^{+}=\left[c_{i}^{+}\left\{c_{j}^{+},c_{k}\right\}\right]
\label{eq57}\\
\left\langle
c_{i},c_{j},c_{k}\right\rangle_{2}&=&0=\left[c_{i}\left\{c_{j},c_{k}\right\}\right]
\label{eq58}\\
\left\langle
c_{i}^{+},c_{j}^{+},c_{k}^{+}\right\rangle_{2}&=&0=\left[c_{i}^{+}\left\{c_{j}^{+},c_{k}^{+}\right\}\right]
\label{eq59}\\
\end{eqnarray}
which reproduce the paraboson statistics.\\
Consider now the case:\\
$A=A_{\left(\begin{array}{c}0
\\
0\end{array}\right)}\oplus A_{\left(\begin{array}{c}1
\\
0\end{array}\right)}\oplus A_{\left(\begin{array}{c}1
\\
1\end{array}\right)}$ and choose the quadratic form $Q$ such that
\begin{eqnarray}
Q(a_{i}^{+},a_{j})&=&Q(a_{j},a_{i}^{+})=\delta_{ij} \nonumber\\
Q(a_{i},a_{j})&=&Q(a_{i}^{+},a_{j}^{+})=0 \nonumber\\
Q(c_{i}^{+},c_{j})&=&-Q(c_{j},c_{i}^{+})=\delta_{ij} \label{eq60}\\
Q(c_{i},c_{j})&=&Q(c_{i}^{+},c_{j}^{+})=0 \nonumber\\
Q(a_{i}^{+},c_{j})&=&Q(a_{j},c_{i}^{+})=Q(a_{i},c_{j})=Q(a_{i}^{+},c_{j}^{+})=0
\nonumber
\end{eqnarray}
$M=\left(\begin{array}{cc} 1 & 0
\\
0 & 1 \end{array}\right)$ like in the precedent examples; we have:
\begin{eqnarray}
\sigma\left(\left(\begin{array}{c}1
\\
0\end{array}\right), \left(\begin{array}{c}1
\\
0\end{array}\right)\right)&=&\left(\begin{array}{cc} 1 &
0\end{array}\right)\left(\begin{array}{cc} 1 & 0
\\
0 & 1 \end{array}\right)\left(\begin{array}{c}1
\\
0\end{array}\right)=1\nonumber\\
\sigma\left(\left(\begin{array}{c}1
\\
1\end{array}\right), \left(\begin{array}{c}1
\\
1\end{array}\right)\right)&=&\left(\begin{array}{cc} 1 &
1\end{array}\right)\left(\begin{array}{cc} 1 & 0
\\
0 & 1 \end{array}\right)\left(\begin{array}{c}1
\\
1\end{array}\right)=2=0\label{eq61}\\
\sigma\left(\left(\begin{array}{c}1
\\
0\end{array}\right), \left(\begin{array}{c}1
\\
1\end{array}\right)\right)&=&\left(\begin{array}{cc} 1 &
0\end{array}\right)\left(\begin{array}{cc} 1 & 0
\\
0 & 1 \end{array}\right)\left(\begin{array}{c}1
\\
1\end{array}\right)=1 \nonumber
\end{eqnarray}
\begin{eqnarray}
\left\langle
a_{i},c_{j},a_{k}^{+}\right\rangle_{1}&=&-2\delta_{ik}c_{j}=\left[\left[a_{i},c_{j}\right],a_{k}^{+}\right]
\label{eq62}\\
\left\langle
a_{i},c_{j},a_{k}^{+}\right\rangle_{1}&=&-2\delta_{ik}c_{j}=\left[a_{i},\left[c_{j},a_{k}^{+}\right]\right]
\label{eq63}
\end{eqnarray}
Note that the bosonic operator $c_{j}$ does not commute with the
fermionic operator $a_{i}$ in Eq. (\ref{eq63}) or with the
fermionic operator $a_{k}^{+}$ in Eq. (\ref{eq64}).\\
The relations below follow easily if .\\
If $A=A_{\left(\begin{array}{c}0
\\
0\end{array}\right)}\oplus A_{\left(\begin{array}{c}1
\\
0\end{array}\right)}\oplus A_{\left(\begin{array}{c}0
\\
1\end{array}\right)}$, $M=\left(\begin{array}{cc} 1 & 0
\\
0 & 1 \end{array}\right)$; let $a\in A_{\left(\begin{array}{c}1
\\
0\end{array}\right)}$, $b\in A_{\left(\begin{array}{c}0
\\
1\end{array}\right)}$ and $c\in A_{\left(\begin{array}{c}0
\\
1\end{array}\right)}$ and
\begin{eqnarray}
Q(a_{i}^{+},a_{j})&=&Q(a_{j},a_{i}^{+})=\delta_{ij} \nonumber\\
Q(a_{i},a_{j})&=&Q(a_{i}^{+},a_{j}^{+})=0 \nonumber\\
Q(b_{i}^{+},b_{j})&=&-Q(b_{j},b_{i}^{+})=\delta_{ij} \label{eq64}\\
Q(b_{i},b_{j})&=&Q(b_{i}^{+},b_{j}^{+})=0 \nonumber\\
Q(a_{i}^{+},b_{j})&=&Q(a_{j},b_{i}^{+})=Q(a_{i},b_{j})=Q(a_{i}^{+},b_{j}^{+})=0
\nonumber
\end{eqnarray}
If we set $a=a_{i}^{+}$, $b=b_{j}^{+}$ and $c=b_{k}$, we find that,
\begin{eqnarray}
\left\langle
a_{i}^{+},b_{j}^{+},b_{k}\right\rangle_{1}&=&2\delta_{jk}a_{i}^{+}=\left[\left[a_{i}^{+},b_{j}^{+}\right],b_{k}\right]
\label{eq65}
\end{eqnarray}
 Since $ A_{\left(\begin{array}{c}1
\\
0\end{array}\right)}$ and $ A_{\left(\begin{array}{c}0
\\
1\end{array}\right)}$ are odd spaces (antisymmetric spaces i.e
fermionic spaces) with respect to $\sigma$, we can
say that different fermion species do not anticommute.\\
If $A=A_{\left(\begin{array}{c}0
\\
0\end{array}\right)}\oplus A_{\left(\begin{array}{c}1
\\
0\end{array}\right)}\oplus A_{\left(\begin{array}{c}0
\\
1\end{array}\right)}$\\
let $a\in A_{\left(\begin{array}{c}1
\\
0\end{array}\right)}$, $b\in A_{\left(\begin{array}{c}0
\\
1\end{array}\right)}$ and $c\in A_{\left(\begin{array}{c}1
\\
0\end{array}\right)}$
\begin{eqnarray}
\left\langle
a,b,c\right\rangle_{1}&=&aQ(b,c)-bQ(a,c)+Q(c,a)b-Q(c,b)a
\nonumber \\
&=&abc-bac+cab-cba \nonumber \\
&=&\left\{\left[a,b\right],c\right\}  \label{eq66}
\end{eqnarray}
we have $Q(b,c)=-Q(c,b)=0$ and $Q(a,c)=Q(c,a)$ then
\begin{eqnarray}
\left\langle
a,b,c\right\rangle_{1}=2Q(b,c)a=\left\{\left[a,b\right],c\right\}
\label{eq67}
\end{eqnarray}
If we set $a=a_{i}^{+}$, $b=b_{j}$ or $b_{j}^{+}$, $c=a_{k}$ and
$Q(a_{i}^{+},a_{k})=\delta_{ik}$ we have
\begin{eqnarray}
\left\langle
a_{i}^{+},b_{j},a_{k}\right\rangle_{1}&=&\left\{\left[a_{i}^{+},b_{j}\right],a_{k}\right\}
\label{eq68} \\
\left\langle
a_{i}^{+},b_{j}^{+},a_{k}\right\rangle_{1}&=&\left\{\left[a_{i}^{+},b_{j}^{+}\right],a_{k}\right\}
\label{eq69}
\end{eqnarray}
\begin{eqnarray}
\left\langle
a,b,c\right\rangle_{2}&=&Q(a,b)c-Q(a,c)b+bQ(c,a)-cQ(b,a)=\left\{a,\left[b,c\right]\right\}
\label{eq70}
\end{eqnarray}
we have $Q(a,b)=-Q(b,a)$ and $Q(a,c)=Q(c,a)$ then
\begin{eqnarray}
\left\langle
a,b,c\right\rangle_{2}&=&2Q(a,b)c=\left\{a,\left[b,c\right]\right\}
\label{eq71}
\end{eqnarray}
If we set $a=a_{i}^{+}$, $b=b_{j}$ or $b_{j}^{+}$, $c=a_{k}$ and
$Q(a_{i}^{+},a_{k})=\delta_{ik}$ we have
\begin{eqnarray}
\left\langle
a_{i}^{+},b_{j},a_{k}\right\rangle_{2}&=&\left\{a_{i}^{+},\left[b_{j},a_{k}\right]\right\}=-2\delta_{ik}b_{j}
\label{eq72} \\
\left\langle
a_{i}^{+},b_{j}^{+},a_{k}\right\rangle_{2}&=&\left\{a_{i}^{+},\left[b_{j}^{+},a_{k}\right]\right\}=-2\delta_{ik}b_{j}^{+}
\label{eq73}
\end{eqnarray}
Note that in this case $\left\langle a,b,c\right\rangle_{1}\neq
\left\langle a,b,c\right\rangle_{2}$\\
If $a\in A_{\left(\begin{array}{c}1
\\
0\end{array}\right)}$, $b\in A_{\left(\begin{array}{c}0
\\
1\end{array}\right)}$ and $c\in A_{\left(\begin{array}{c}1
\\
1\end{array}\right)}$
\begin{eqnarray}
\left\langle
a,b,c\right\rangle_{1}&=&aQ(b,c)+bQ(a,c)-Q(c,a)b-Q(c,b)a
\nonumber \\
&=&abc+bac-cab-cba \nonumber \\
&=&\left\{\left[a,b\right],c\right\}  \label{eq74}
\end{eqnarray}
we have
\begin{eqnarray}
Q(b,c)=Q(c,b)~~~~~~~~\mbox{and}~~~~~~~~Q(a,c)=Q(c,a)~~~~~~~~\mbox{
then}\label{eq75}
\end{eqnarray}
\begin{eqnarray}
\left\langle
a,b,c\right\rangle_{1}=0=\left\{\left[a,b\right],c\right\}
\label{eq76}
\end{eqnarray}
\begin{eqnarray}
\left\langle
a,b,c\right\rangle_{2}&=&Q(a,b)c-Q(a,c)b+bQ(c,a)b-cQ(b,a)
\nonumber \\
&=&abc-acb+bca-cba \nonumber \\
&=&\left\{a,\left[b,c\right]\right\}  \label{eq77}
\end{eqnarray}
we have
\begin{eqnarray}
Q(a,b)=-Q(b,a)~~~~~~~~\mbox{and}~~~~~~~~Q(a,c)=Q(c,a)~~~~~~~~\mbox{
then}\label{eq78}
\end{eqnarray}
\begin{eqnarray}
\left\langle
a,b,c\right\rangle_{2}=2Q(a,b)c=\left\{a,\left[b,c\right]\right\}
\label{eq79}
\end{eqnarray}
Note that $a$, $b$ and $c$ belong to different spaces, by definition
$Q(a,b)=0$, $Q(a,c)=0$ and $Q(b,c)=0$, hence:
\begin{eqnarray}
\left\langle
a,b,c\right\rangle_{2}=2Q(a,b)c=\left\{a,\left[b,c\right]\right\}=0
\label{eq80}
\end{eqnarray}
%Let
%\begin{eqnarray}
%M&=&\left(\begin{array}{cc} \alpha_{11} & \alpha_{12}
%\\
%\alpha_{21} & \alpha_{22}
%\end{array}\right)=\left(\begin{array}{cc} 0 &
%1
%\\
%1 & 0 \end{array}\right)\label{eq81}\\
%\sigma\left(\left(\begin{array}{c}1
%\\
%0\end{array}\right), \left(\begin{array}{c}1
%\\
%0\end{array}\right)\right)&=&\left(\begin{array}{cc} 1 &
%0\end{array}\right)\left(\begin{array}{cc} 0 & 1
%\\
%1 & 0 \end{array}\right)\left(\begin{array}{c}1
%\\
%0\end{array}\right)=0\label{eq82}\\
%\forall x, y\in A_{\left(\begin{array}{c}1
%\\
%0\end{array}\right)}&,& Q(x,y)=-Q(y,x) \label{eq83}
%\end{eqnarray}

\section{Application to the bilinear case}
\subsection{Anyonic Vector Spaces}
We consider braided categories associated to $\mathbb Z_n$, the finite group of order n. Let $g$ be the generator of $\mathbb Z_n$ with $g^{n}~=~1$. As a category of objects and morphisms we take the category $\bf{Rep}(\mathbb Z_n)$ of finite dimensional representations of $\mathbb Z_n$. \\
Given an object $V$ of $\bf{Rep}(\mathbb Z_n)$ we can decompose it
under the action of $\mathbb Z_n$ as
\begin{eqnarray}
V~=~\oplus_{p=0}^{n-1}V_p,~~~~~~~~~~~p~=~0,1\ldots,n-1.
\end{eqnarray}
Here $p$ runs over the set of irreducible representation $\rho_p$. We have the action,\\
$~~~~~~~~~~~~~~~~~~~~~~~~~~~~g\triangleright v~=~e^{\frac{2\pi i p}{n}} v,~~~~~~~~~~~~~~~~\forall~ v \in V_p$\\
where the action of $\mathbb Z_n$ is denoted by~$ \triangleright$. If$ v~\in V_p$, we say that $v$ is homogeneous of degree $|v|~=~p$\\
On this category $\bf{Rep}(\mathbb Z_n)$, we can now define the
braiding
\begin{eqnarray}
\Psi_{V,W}( v \otimes w )~=~e^{\frac{2\pi i |v||w|}{n}}w \otimes v
\end{eqnarray}
In physics the quantities~ $e^{\frac{2\pi i |v||w|}{n}}$~are called
fractional or anyonic statistics.
\subsection{Classical case}
$A=\oplus_{i}A_{v_{i}}$ is a $\Gamma$ -graded algebra$^{5}$, we
consider (anti)symmetric bilinear maps that we have defined in an
earlier paper$^{5}$ $\sigma :\Gamma\times\Gamma\longrightarrow
Z_{2}; Z_{2}=\left\{\overline{0},\overline{1}\right\}$; and $T(A)$
the tensor algebra of $A$. We define isomorphisms
$\Psi_{A_{v_{i}},A_{v_{j}}} :A_{v_{i}}\otimes
A_{v_{j}}\longrightarrow A_{v_{j}}\otimes A_{v_{i}}$ which satisfy
the properties (2). Let $\left[~,~\right]_{\sigma}$ be a bilinear
map, $\left[~,~\right]_{\sigma} :A\times A\longrightarrow A$; and
$J$ be the two sided ideal generated by all expressions of the
following form, $x\otimes y+\Psi_{A_{v_{i}},A_{v_{j}}}(x\otimes
y)-\left[x,y\right]_{\sigma}$ where $x\in A_{v_{i}}$, $y\in
A_{v_{j}}$, $i$ and $j\in Z_{2}$.
The envelopping in algebra is defined as $U(A)=T(A)/J$.\\
We take the braiding as:\\
\begin{eqnarray}
\Psi_{A_{v_{i}},A_{v_{j}}}(x\otimes
y)=-(-1)^{\sigma(v_{i},v_{j})}y\otimes x \label{eqq26}
\end{eqnarray}
In $U(A)$, and also in $A$ we have:
\begin{eqnarray}
\left[x,y\right]_{\sigma}=x.y-(-1)^{\sigma(v_{i},v_{j})}y.x
\label{eqq27}
\end{eqnarray}
this is the case that we have developped in the paper$^{5}$, which
gives the color superalgebra, the Lie superalgebra, the color
algebra, etc...\\
Note that the anyonic case is obtained by taking the braiding as:
\begin{eqnarray}
\Psi_{A_{v_{i}},A_{v_{j}}}(x\otimes y)=exp({\frac{2\pi
i\sigma(v_{i},v_{j})}{n}})y\otimes x \label{eqq28}
\end{eqnarray}
now the bilinear map $\sigma$ is:
\begin{eqnarray}
\sigma :\Gamma\times\Gamma\longrightarrow Z_{n};
Z_{n}=\left\{\overline{0},\overline{1},...,\overline{n-1}\right\}\label{eqq29}
\end{eqnarray}
Then in $U(A)$, and also in $A$, we have:
\begin{eqnarray}
\left[x,y\right]_{\sigma}=x.y+exp({\frac{2\pi
i\sigma(v_{i},v_{j})}{n}})y.x \label{eqq30}
\end{eqnarray}
\section{Schur functors}
It is well known that the functors generalizing the symmetric powers
an exterior powers are defined in terms of the Young symmetrizers
$c_\lambda$. For any finite complex vector space $A$,  %given representation $A$ of an arbitrary group $G$,
we consider the dth tensor power of $A$, on which the symmetric
group $\mathcal{S}_d$ acts, say on the right, by permuting the
factors $(v_1\otimes ...\otimes
v_d).\sigma=v_{\sigma(1)}\otimes...\otimes v_{\sigma(d)}$. This
action commute with the left action of $GL(A)$. For any partition
$\lambda$ of $d$ we have a Young Symmetrizer $c_\lambda$ in the
group algebra $\mathbb{C}\mathcal{S}_d$. We denote the image of
$c_\lambda$ on $A^{\otimes d}$ by $\mathbb{S}_{\lambda}A$ . \\( The
functoriality means that a linear map $\varphi : V ~\rightarrow~W$
of vector spaces determines a linear map
~$\mathbb{S}_{\lambda}(\varphi):\mathbb{S}_{\lambda}(V)~\rightarrow~\mathbb{S}_{\lambda}(W)$,
with
$\mathbb{S}_{\lambda}(\varphi\circ\psi)~=~\mathbb{S}_{\lambda}(\varphi)\circ~\mathbb{S}_{\lambda}(\psi)$
and $\mathbb{S}_{\lambda}(Id_{V})~=~Id_{\mathbb{S}_{\lambda} V}$)

\begin{eqnarray}
\mathbb{S}_{\lambda}A~=~Im(~c_{\lambda}\vert_ {A^{\otimes d}})
\end{eqnarray}
Which is a representation of $GL(A)$. The functor
$A~\rightsquigarrow~\mathbb{S}_{\lambda}A$ is called the Schur
functor or Weyl module corresponding to $\lambda$. We have the
canonical decomposition
\begin{eqnarray}
A\otimes A ~=~Sym^{2}A~\oplus~\Lambda^{2}A
\end{eqnarray}
The group $GL(A)$ acts on $A\otimes A$, and decompose it into a
direct sum of irreducible $GL(A)$-representations. We will be
interested by the next tensor power
\begin{eqnarray}
A\otimes A\otimes A ~=~Sym^{3}A~\oplus~\Lambda^{3}A
\oplus~(\mathbb{S}_{(2,1)}A)^{\oplus 2}
\end{eqnarray}
The partition d~=~3 corresponds to the functor
$A~\rightsquigarrow~Sym^{3}A$ and the partition\\$d~=~1+1+1$ to the
functor $A~\rightsquigarrow~\Lambda^{3}A$. We found new and more
generalized things; for example the partition 3~=~2+1, the
corresponding  "symmetrizer" $c_{\lambda,\Psi}$ is
\begin{eqnarray}
c_{(2,1),\Psi}~=~e_1 +~
\Psi_{ij}^{mn}~e_{12}+\Psi_{ij}^{mn}\Psi_{nk}^{pq}\Psi_{mp}^{rs}~e_{13}+\Psi_{ij}^{mn}\Psi_{nk}^{pq}\Psi_{mp}^{rs}\Psi_{sq}^{tu}~e_{132}
\end{eqnarray}
\begin{eqnarray}
c^{*}_{(2,1),\Psi}=e_1+~\Psi_{jk}^{mn}e_{23}
+\Psi_{jk}^{mn}\Psi_{im}^{pq}\Psi_{qn}^{rs}~e_{13}+\Psi_{jk}^{mn}\Psi_{im}^{pq}\Psi_{qn}^{rs}\Psi_{pr}^{tu}~e_{123}
\end{eqnarray}
\section{Conclusion}
We can consider the algebra $A\otimes \mathbb{C}^\infty
(\mathbb{R}^n)$ and interested by $S_3$ irreducible representations
of such algebra ; this algebra may be viewed as an algebra of
operators, and also to understand the unitarity condition (\ref{e6}) in this case.\\
\textbf{Acknowledgment}:\\
I thank the professor A. Elmazouni for helpful mathematical
discussions.


\begin{thebibliography}{99}
\bibitem{gr}
{\sc H. S. Green},
 Phys. Rev. 90, 270(1953).
\bibitem{ne}
{\sc A. Neveu and J. H. Schwarz, Nucl. Phy. B 31, 86 (1971)}.
\bibitem{Ka}
{\sc V.G. Kac, Adv.Math. 26, 8 (1977).}
\bibitem{Oh}
{\sc Y.Ohnuki and S.Kamafuchi, Quantum Field Theory and
Parastatistics (University of Tokyo) Tokyo(1982)}.
\bibitem{Zi}
{\sc A.Zinoun and J.Cortois, J.Math.Phys. 32, 247 (1991)}.
\bibitem{Bi}
{\sc S.N.Biswa and S.K.Soni, J.Math.Phys. 29, 16 (1988)}.
\bibitem{Ma}
{\sc S.Majid, J.Math.Phys. 32, 12 (1991)}.
\bibitem{Ca}
{\sc C.Cassel, Publication de l'institut de recherche Mathématique
avancée. R.C.P 25, vol 43, Strasbourg}.
\bibitem{PC}
{\sc P.Cartier, An introduction to quantum groups. IHES/M/7(1992)}.

\end{thebibliography}
\end{document}